\documentclass[11pt,fleqn]{elsarticle}
\usepackage{amssymb}

\usepackage{hyperref}
\usepackage{times}
\usepackage{graphicx}
\usepackage{wrapfig}
\usepackage{fancyhdr}
\usepackage{amsmath}

\input{tcilatex}
\begin{document}

\begin{frontmatter}

\title{Laplace Green's Functions for Infinite Ground Planes with Local Roughness}
	\author{Nail A. Gumerov and Ramani Duraiswami \\
		University of Maryland Institute for Advanced Computer Studies, College Park, MD 20742}

		\begin{keyword}
			Green's Function\sep Laplace's Equation \sep Electrostatics \sep infinite ground
			\MSC[2010] 00-01\sep  99-00
		\end{keyword}

\begin{abstract}
	The Green's functions for the Laplace equation respectively satisfying the Dirichlet and Neumann boundary conditions on the upper side of an infinite plane with a circular hole are introduced and constructed. These functions enables solution of the boundary value problems in domains where the hole is closed by any surface. This approach enables accounting for arbitrary positive and negative ground elevations inside the domain of interest, which, generally, is not possible to achieve using the regular method of images. Such problems appear in electrostatics, however, the methods developed apply to other domains where the Laplace or Poisson equations govern. Integral and series representations of the Green's functions are provided. An efficient computational technique based on the boundary element method with fast multipole acceleration is developed. A numerical study of some benchmark problems is presented.
\end{abstract}
	\end{frontmatter}


\section{Introduction}

Many problems of physical interest involve fields satisfying Laplace's
equation with appropriate boundary conditions. Often, the boundaries of the
domain might include surfaces of infinite extent, which for analytical
convenience, are often considered to be infinite planes, so that they can be
treated via the method of images. Examples include electrostatics,
magnetostatics, conduction (or diffusion) from a surface at a fixed
temperature (or concentration), potential flow of an ideal fluid over a
body. Motivating the work reported here are the problems of electric and
magnetic fields near an infinite ground, where the surface elevation departs
from the plane in both the positive and negative directions, e.g., a ground
with bumps and troughs.

Computation of electrostatic fields in scales where the effects of ground
should be taken into account can be an important task for a number of
applications, for example, for modeling fields in urban environments \cite%
{Adelman2017:SC}, fields generated by the transmission power lines \cite%
{DAmore96:IEEE}, \cite{Trlep2009:IEEE}, and fields generated by lightning
over rough ocean surface \cite{Zhang2012:Electrostat}. Scattering of higher
frequency electromagnetic waves from rough surfaces is also of interest \cite%
{Chen2003:IEEE}. Note that if the wavelength is much larger than the scale
of the roughness, electrostatic approximation for characterization of
scattering properties of the surface is also valid.

Most of numerical methods used for such computations either neglect the
presence of objects outside the computational domain or impose some boundary
conditions on the domain surface (e.g., see monograph \cite{Givoli1992:book}
and review \cite{Tsinkov1998:ANM}). The boundary elements can deliver
solutions which provide decay of the electric potential at the infinity,
while solving a boundary integral equation for the charge distribution only
for the objects located inside the computational domain. It is relatively
easy to account the effects of the infinite flat ground in such computations
using the method of images \cite{Jackson1998:Wiley}. Formally, this can be
done by replacing the free space Green's function in the boundary integral
equation by the Green's function for the semispace. This function takes zero
at any point on the ground plane and so any convolution with such a Green's
function produces a solution satisfying the boundary condition on the
ground. Some techniques for construction of Green's functions for the
Laplace equation in different 2D domains are developed (e.g., \cite%
{Crowdy2007:IMA}), but they formulated in terms of complex variables and
conformal mappings, which makes difficult to extend them for 3D.

In practice, the situation frequently appears to be more complicated as the
ground has elevations and dips. The term ``ground'' can also include
buildings or other construction, and can be modeled using surface boundary
element meshes. In any case, such modeling is limited to the area where the
solution is needed and detailed representation of the environment can be
provided only inside the computational domain. The question then is how to
take into account the effect of the environment outside the computational
domain? One possibility is to assume periodic boundary conditions. In this
case some ``periodization'' techniques for solution of the Laplace equation
can be used (e.g., see \cite{Gumerov2014JCP}; quasi-periodic Green's
functions can be found in \cite{Moroz2006:JP}). Another approach can be to
model the environment outside some ball $\mathbf{B}$ of radius $R$ as a flat
ground, while provide a detailed meshing of all objects and the ground
inside $\mathbf{B}$ (sometimes called ``locally rough surface''). In the
present paper we take the latter approach.

A simple solution of the problem can be found if all points on the ground
inside $\mathbf{B}$ are located above or on the level of the flat ground at
infinity (which we can denote as zero level). In this case the method of
images can be applied to any objects in the domain. However, if there exist
some dips, or points below the zero level the method of images is not
applicable, since it requires placing the image source on the side of the
ground plane opposite to the side where the actual source is located, which
creates a physically unacceptable (singular) solution. In a number of papers
treatment of such situation using boundary integrals for the Helmholtz
equation was considered (e.g., \cite{Bao2013:IPI}, \cite%
{Chandler-Wilde2006:SIAM}). In this paper, we take a different approach and
derive an analogue of the Green's function, which accounts for the effect of
infinite flat ground plane outside the computational domain, and which has
no singularities besides the source location inside $\mathbf{B}$. This
function is suitable for solution of problems with arbitrary location of
boundaries inside $\mathbf{B}$. We implemented the method and provide some
numerical examples. Note that this method can also be applied to the
Helmholtz equation.

The problem described above formally reduces to solution of the Laplace
equation with the Dirichlet boundary conditions on the ground. In practice,
there appear also situations when the zero Neumann boundary conditions
should be imposed on an infinite surface. For example, if we wish to
determine underwater electric fields generated by some sources, this can be
formulated as the problem with zero Neumann boundary conditions on the sea
level (since the electric conductivity of the sea water is much larger than
that of the air). The ocean surface can be rough, so that can be taken into
account for points close to the surface, while the ocean surface can be
considered as flat in the far field (also, see \cite{Holmes2015:MTSJ}, \cite%
{Yue2016:IEEE}, \cite{Wang2018:PIER}). Another example of the Neumann
problem appears in aerodynamics, where the effects of ground should be taken
into account to determine correctly the lift force for objects flying close
to the ground \cite{Rozhdestvensky2000:book}. For flat surfaces the Neumann
problem can be also solved using the method of images, which causes the same
problems as for the Dirichlet problem described above. The Green's functions
for the Dirichlet and Neumann boundary conditions are closely related. In
this paper we address both problems.

\section{Statement of the problem}

\begin{figure}[tbh]
\vspace{-20pt}
\par
\begin{center}
\includegraphics[width=0.9\textwidth, trim=0 1.15in 1.5in 0]{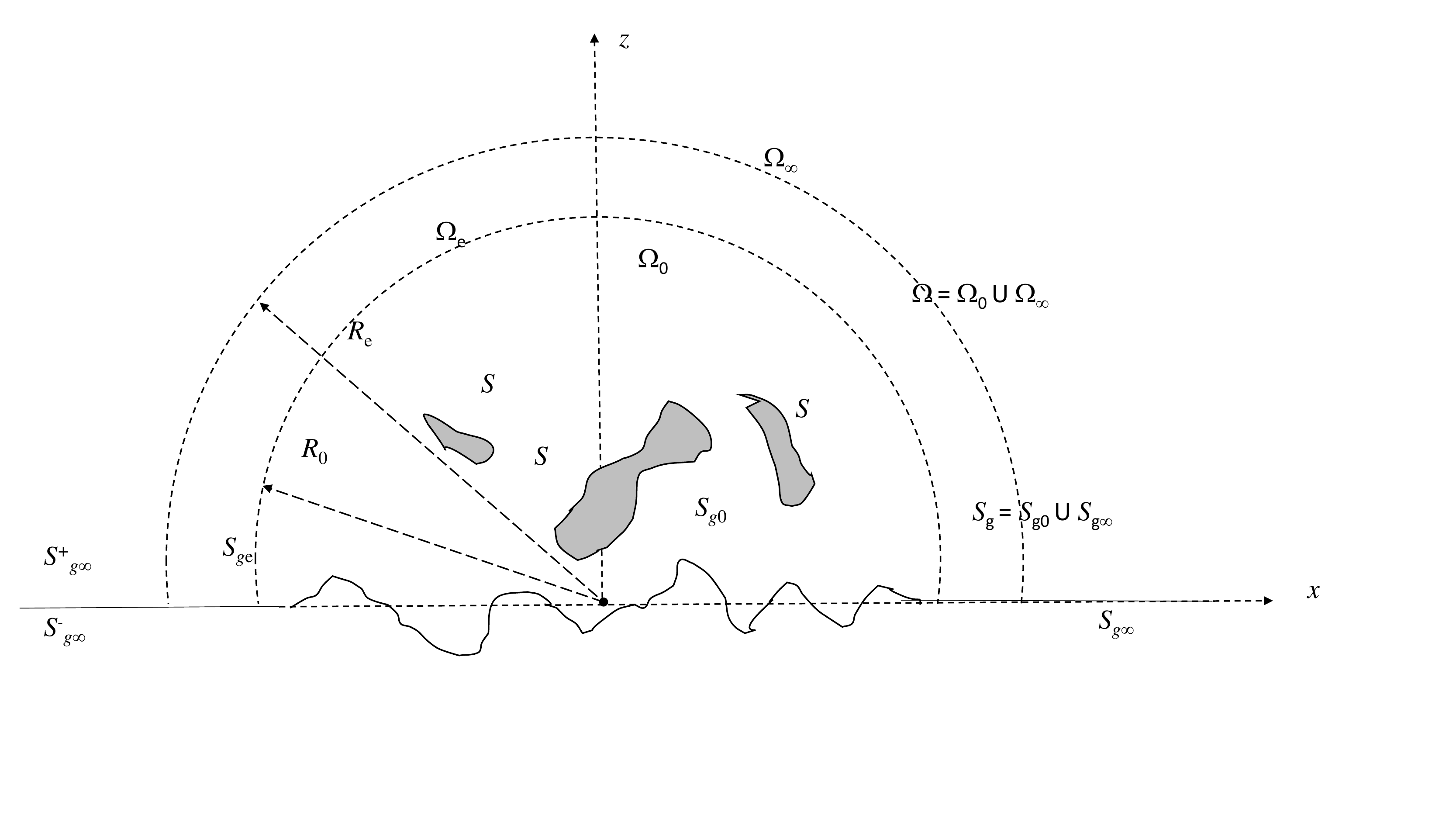}
\end{center}
\caption{The problem geometry and notation. A solution needs to be obtained
in domain $\Omega _{0}$ bounded by the object surface $S$, ground surface $%
S_{g0}$ and a hemisphere of radius $R_{0}$. The ground surface is flat
outside $\Omega _{0}$. An extended computational domain $\Omega _{e}$, which
is a part of $\Omega $ located inside a sphere of radius $R_{e}$ can be
introduced for efficient computations. }
\label{Fig1}
\end{figure}

First, we consider the Dirichlet problem. Consider an infinite ground
surface $S_{g}$ on which the electric potential $\phi $ is zero and let
there be some objects on it. Let the surface be $S$ the value $V$ of the
potential is given (here and below we use dimensionless equations). In the
electrostatics approximation, we have the following Dirichlet problem for
the Laplace equation, 
\begin{eqnarray}
\nabla ^{2}\phi \left( \mathbf{y}\right) &=&0,\quad \mathbf{y}\in \Omega ,
\label{sp1} \\
\left. \phi \right| _{\mathbf{y\in }S_{g}} &=&0,\quad \left. \phi \right| _{%
\mathbf{y\in }S}=V,\quad \left. \phi \right| _{\left| \mathbf{y}\right| 
\mathbf{=}\infty }=0,  \notag
\end{eqnarray}%
where $\Omega $ is an infinite domain in $\mathbb{R}^{3}$ bounded by $S$ and 
$S_{g}$. The solution of this problem can be obtained using indirect
boundary integral method, which represents the potential in the form of the
single layer potential 
\begin{equation}
\phi \left( \mathbf{y}\right) =\int_{S\cup S_{g}}\sigma \left( \mathbf{x}%
\right) G\left( \mathbf{y},\mathbf{x}\right) dS\left( \mathbf{x}\right)
,\quad G\left( \mathbf{y},\mathbf{x}\right) =\frac{1}{4\pi \left| \mathbf{y}-%
\mathbf{x}\right| },  \label{sp2}
\end{equation}%
where $G\left( \mathbf{y},\mathbf{x}\right) $ is the free-space Green's
function for the Laplace equation, and $\sigma \left( \mathbf{x}\right) $ is
the surface charge density, which can be obtained from the solution of the
boundary integral equation, 
\begin{equation}
\int_{S\cup S_{g}}\sigma \left( \mathbf{x}\right) G\left( \mathbf{y},\mathbf{%
x}\right) dS\left( \mathbf{x}\right) =\left\{ 
\begin{array}{c}
V\left( \mathbf{y}\right) ,\quad \mathbf{y}\in S, \\ 
0,\quad \mathbf{y}\in S_{g}.%
\end{array}%
\right.  \label{sp3}
\end{equation}%
Assume further that $S_{g}$ can be modelled by a flat horizontal surface
(zero elevation, $z=0$), except for a finite part near the source, and in
the region in which we wish to compute the solution in details. The ground
plane in this region may have arbitrary positive and negative elevations. We
assume that this region is contained inside a ball $\mathbf{B}_{0}$ of
radius $R_{0}$ centered at the origin of the reference frame. It is also
assumed that surface $S$ is located inside the ball. The ball partitions $%
\Omega $ into domains $\Omega _{0}$ and $\Omega _{\infty }$ and $S_{g}$ into 
$S_{g0}$ and $S_{g\infty }$ (see Fig. \ref{Fig1}), 
\begin{eqnarray}
\Omega _{0} &=&\left\{ \mathbf{y}\in \Omega ,\quad \left| \mathbf{y}\right|
<R_{0}\right\} ,\quad \Omega _{\infty }=\left\{ \mathbf{y}\in \Omega ,\quad
\left| \mathbf{y}\right| \geqslant R_{0}\right\} =\left\{ \mathbf{y},\quad
\left| \mathbf{y}\right| \geqslant R_{0},\text{ }z>0\right\} ,  \label{sp3.1}
\\
S_{g0} &=&\left\{ \mathbf{y}\in S_{g},\quad \left| \mathbf{y}\right|
<R_{0}\right\} ,\quad S_{g\infty }=\left\{ \mathbf{y}\in S_{g},\quad \left| 
\mathbf{y}\right| \geqslant R_{0}\right\} =\left\{ \mathbf{y},\quad \left| 
\mathbf{y}\right| \geqslant R_{0},\text{ }z=0\right\} .\mathbf{\quad } 
\notag
\end{eqnarray}%
Note now that $S_{g\infty }$ is an open surface, for which we have an
``upper'' side, $S_{g\infty }^{+}$, faced towards $z>0$, and a ``lower''
side, $S_{g\infty }^{-}$ faced towards $z<0$. The boundary conditions need
be satisfied on the upper side only. The Green's function $G^{(D)}\left( 
\mathbf{y},\mathbf{x}\right) $ for the zero Dirichlet boundary conditions on 
$S_{g\infty }^{+}$ can be introduced as the solution of the problem 
\begin{eqnarray}
\nabla _{\mathbf{y}}^{2}G^{(D)}\left( \mathbf{y},\mathbf{x}\right)
&=&-\delta \left( \mathbf{y-x}\right) ,\quad \quad \mathbf{y,x}\in \mathbb{R}%
^{3}\mathbf{,}  \label{sp4} \\
\left. G^{(D)}\left( \mathbf{y},\mathbf{x}\right) \right| _{\mathbf{y\in }%
S_{g\infty }^{+},\mathbf{y}\neq \mathbf{x}} &=&0,\text{\quad }\left.
G^{(D)}\left( \mathbf{y},\mathbf{x}\right) \right| _{\left| \mathbf{y}%
\right| \mathbf{=}\infty ,\mathbf{y}\neq \mathbf{x}}=0.  \notag
\end{eqnarray}%
The solution of the Dirichlet problem (\ref{sp1}) then can be sought in the
form 
\begin{equation}
\phi \left( \mathbf{y}\right) =\int_{S\cup S_{g0}}\sigma ^{(D)}\left( 
\mathbf{x}\right) G^{(D)}\left( \mathbf{y},\mathbf{x}\right) dS\left( 
\mathbf{x}\right) ,  \label{sp5}
\end{equation}%
where $\sigma ^{(D)}\left( \mathbf{x}\right) $ satisfies the boundary
integral equation 
\begin{equation}
\int_{S\cup S_{g0}}\sigma ^{(D)}\left( \mathbf{x}\right) G^{(D)}\left( 
\mathbf{y},\mathbf{x}\right) dS\left( \mathbf{x}\right) =\left\{ 
\begin{array}{c}
V\left( \mathbf{y}\right) ,\quad \mathbf{y}\in S, \\ 
0,\quad \mathbf{y}\in S_{g0}.%
\end{array}%
\right.  \label{sp6}
\end{equation}%
We seek to derive analytical expressions for, and methods for computation
of, the function $G^{(D)}\left( \mathbf{y},\mathbf{x}\right) $.

Similarly, we can also get the Green's function $G^{(N)}\left( \mathbf{y},%
\mathbf{x}\right) $ for the zero Neumann boundary conditions on $S_{g\infty
} $, which can be found by solving 
\begin{eqnarray}
\nabla _{\mathbf{y}}^{2}G^{(N)}\left( \mathbf{y},\mathbf{x}\right)
&=&-\delta \left( \mathbf{y-x}\right) ,\quad \quad \mathbf{y,x}\in \mathbb{R}%
^{3}\mathbf{,}\quad  \label{sp7} \\
\left. \frac{\partial }{\partial n\left( \mathbf{y}\right) }G^{(N)}\left( 
\mathbf{y},\mathbf{x}\right) \right| _{\mathbf{y\in }S_{g\infty }^{+},%
\mathbf{y}\neq \mathbf{x}} &=&0,\text{\quad }\left. G^{(N)}\left( \mathbf{y},%
\mathbf{x}\right) \right| _{\left| \mathbf{y}\right| \mathbf{=}\infty , 
\mathbf{y}\neq \mathbf{x}}=0.  \notag
\end{eqnarray}
Construction of this function allows one to solve the problem with mixed
boundary conditions 
\begin{eqnarray}
\nabla ^{2}\phi \left( \mathbf{y}\right) &=&0,\quad \mathbf{y}\in \Omega ,
\label{sp8} \\
\left. \frac{\partial }{\partial n\left( \mathbf{y}\right) }\phi \right| _{%
\mathbf{y\in }S_{g}} &=&0,\quad \left. \phi \right| _{\mathbf{y\in }%
S}=V,\quad \left. \phi \right| _{\left| \mathbf{y}\right| \mathbf{=}\infty
}=0,  \notag
\end{eqnarray}
in form (\ref{sp5}), where superscript $\left( D\right) $ should be changed
to $\left( N\right) $ and the charge density can be determined from the
boundary integral equation 
\begin{eqnarray}
\int_{S\cup S_{g0}}\sigma ^{(N)}\left( \mathbf{x}\right) G^{(N)}\left( 
\mathbf{y},\mathbf{x}\right) dS\left( \mathbf{x}\right) &=&V\left( \mathbf{y}%
\right) ,\quad \mathbf{y}\in S,  \label{sp9} \\
\frac{\partial }{\partial n\left( \mathbf{y}\right) }\int_{S\cup
S_{g0}}\sigma ^{(N)}\left( \mathbf{x}\right) G^{(N)}\left( \mathbf{y}, 
\mathbf{x}\right) dS\left( \mathbf{x}\right) &=&0,\quad \mathbf{y}\in S_{g0}.
\notag
\end{eqnarray}

Note that functions $G^{(D)}\left( \mathbf{y},\mathbf{x}\right) $ and $%
G^{(N)}\left( \mathbf{y},\mathbf{x}\right) $ defined above are not unique,
since the boundary condition is specified only on $S_{g\infty }^{+}$, while
an arbitrary boundary condition can be imposed on $S_{g\infty }^{-}$.
Depending on the type of condition on $S_{g\infty }^{-}$ the functions may
or may not satisfy the symmetry typical for Green's functions for Dirichlet
problems for closed or two-sided open surfaces with homogeneous boundary
conditions \cite{Duffy2015book}. Despite this non-uniqueness \emph{any}
solution of problems (\ref{sp4}) and (\ref{sp7}) delivers the Green's
function suitable for representation of solutions in form (\ref{sp5} ),
where the charge density is determined from Eqs (\ref{sp6}) and (\ref{sp9}).
Indeed, these relations guarantee that the solution satisfies the Laplace
equation and the boundary conditions, and, so, as $G^{(D)}$ or $G^{(N)}$ are
selected the solution is unique.

We also remark that, despite $G^{(D)}\left( \mathbf{y},\mathbf{x}\right) $
and $G^{(N)}\left( \mathbf{y},\mathbf{x}\right) $ can be computed for any
location of the source ($\mathbf{x}$) and evaluation ($\mathbf{y}$) points,
for boundary integral representation of the solution we only consider the
case when the source $\mathbf{x}$ is located inside the ball $\mathbf{B}$.
Moreover, $\mathbf{y}$ should be also located in $\mathbf{B}$ to solve the
boundary integral equations (\ref{sp6}) and (\ref{sp9}). As soon as the
charge density is determined from these equations the values of Green's
functions at arbitrary $\mathbf{y}$ can be used.

\section{Solution}

\subsection{Single and double layer potentials}

Solutions of the Dirichlet and Neumann problems (\ref{sp4}) and (\ref{sp7})
can be found using single and double layer potentials, 
\begin{eqnarray}
G^{(D)}\left( \mathbf{y},\mathbf{x}\right) &=&G\left( \mathbf{y},\mathbf{x}%
\right) +K^{(D)}\left( \mathbf{y},\mathbf{x}\right) ,\quad K^{(D)}\left( 
\mathbf{y,x}\right) =\int_{S_{g\infty }}\mu _{D}\left( \mathbf{x}^{\prime }%
\mathbf{,x}\right) \frac{\partial G\left( \mathbf{y},\mathbf{x}^{\prime
}\right) }{\partial n\left( \mathbf{x}^{\prime }\right) }dS\left( \mathbf{x}%
^{\prime }\right) ,  \label{so1} \\
G^{(N)}\left( \mathbf{y},\mathbf{x}\right) &=&G\left( \mathbf{y},\mathbf{x}%
\right) +K^{(N)}\left( \mathbf{y},\mathbf{x}\right) ,\quad K^{(N)}\left( 
\mathbf{y,x}\right) =\int_{S_{g\infty }}\sigma _{N}\left( \mathbf{x}^{\prime
}\mathbf{,x}\right) G\left( \mathbf{y},\mathbf{x}^{\prime }\right) dS\left( 
\mathbf{x}^{\prime }\right) .\quad  \notag
\end{eqnarray}%
Here $\sigma _{N}\left( \mathbf{x}^{\prime }\mathbf{,x}\right) $ and $\mu
_{D}\left( \mathbf{x}^{\prime }\mathbf{,x}\right) $ are the single and
double layer densities, which should be found from the boundary conditions.
These conditions show that 
\begin{eqnarray}
\left. K^{(D)}\left( \mathbf{y},\mathbf{x}\right) \right| _{\mathbf{y\in }%
S_{g\infty }^{+}} &=&-\left. G\left( \mathbf{y},\mathbf{x}\right) \right| _{%
\mathbf{y\in }S_{g\infty }},  \label{so1.1} \\
\left. \frac{\partial K^{(N)}\left( \mathbf{y},\mathbf{x}\right) }{\partial
n\left( \mathbf{y}\right) }\right| _{\mathbf{y\in }S_{g\infty }^{+}}
&=&-\left. \frac{\partial G\left( \mathbf{y},\mathbf{x}\right) }{\partial
n\left( \mathbf{y}\right) }\right| _{\mathbf{y\in }S_{g\infty }}.  \notag
\end{eqnarray}

The single and double layer potentials (\ref{so1}) obey the following
symmetry property 
\begin{equation}
\left. K^{(D)}\left( \mathbf{y,x}\right) \right| _{\mathbf{y\in }S_{g\infty
}^{-}}=-\left. K^{(D)}\left( \mathbf{y,x}\right) \right| _{\mathbf{y\in }
S_{g\infty }^{+}},\quad \left. K^{(N)}\left( \mathbf{y,x}\right) \right| _{%
\mathbf{y\in }S_{g\infty }^{-}}=\left. K^{(N)}\left( \mathbf{y,x}\right)
\right| _{\mathbf{y\in }S_{g\infty }^{+}}.  \label{so2}
\end{equation}
Furthermore, we note that the single and double layer densities are related
to the jumps of the potentials and their normal derivatives on the surface
as 
\begin{eqnarray}
\mu _{D}\left( \mathbf{y,x}\right) &=&\left. K^{(D)}\left( \mathbf{y,x}
\right) \right| _{\mathbf{y\in }S_{g\infty }^{+}}-\left. K^{(D)}\left( 
\mathbf{y,x}\right) \right| _{\mathbf{y\in }S_{g\infty }^{-}},  \label{so3}
\\
\sigma _{N}\left( \mathbf{y},\mathbf{x}\right) &=&\left. \frac{\partial
K^{(N)}\left( \mathbf{y},\mathbf{x}\right) }{\partial n\left( \mathbf{y}
\right) }\right| _{\mathbf{y\in }S_{g\infty }^{-}}-\left. \frac{\partial
K^{(N)}\left( \mathbf{y},\mathbf{x}\right) }{\partial n\left( \mathbf{y}
\right) }\right| _{\mathbf{y\in }S_{g\infty }^{+}}.  \notag
\end{eqnarray}
Symmetry (\ref{so2}) then provides 
\begin{eqnarray}
\mu _{D}\left( \mathbf{y,x}\right) &=&2\left. K^{(D)}\left( \mathbf{y}, 
\mathbf{x}\right) \right| _{\mathbf{y\in }S_{g\infty }^{+}},  \label{so4} \\
\sigma _{N}\left( \mathbf{y},\mathbf{x}\right) &=&-2\left. \frac{\partial
K^{(N)}\left( \mathbf{y},\mathbf{x}\right) }{\partial n\left( \mathbf{y}
\right) }\right| _{\mathbf{y\in }S_{g\infty }^{+}}.  \notag
\end{eqnarray}%
Boundary conditions (\ref{so1.1}) complete the solution, which can be
written in the integral form 
\begin{eqnarray}
K^{(D)}\left( \mathbf{y,x}\right) &=&-2\int_{S_{g\infty }}G\left( \mathbf{x}%
^{\prime }\mathbf{,x}\right) \frac{\partial G\left( \mathbf{y},\mathbf{x}%
^{\prime }\right) }{\partial n\left( \mathbf{x}^{\prime }\right) }dS\left( 
\mathbf{x}^{\prime }\right) ,  \label{so5} \\
K^{(N)}\left( \mathbf{y,x}\right) &=&2\int_{S_{g\infty }}G\left( \mathbf{y},%
\mathbf{x}^{\prime }\right) \frac{\partial G\left( \mathbf{x}^{\prime }, 
\mathbf{x}\right) }{\partial n\left( \mathbf{x}^{\prime }\right) }dS\left( 
\mathbf{x}^{\prime }\right) .\quad  \notag
\end{eqnarray}

These equations show that $\mathbf{x}$ and $\mathbf{y}$ enter the expression
in a non-symmetric way. However, \ due to the symmetry of the free-space
Green's function, we have 
\begin{equation}
K^{(N)}\left( \mathbf{y,x}\right) =-K^{(D)}\left( \mathbf{x,y}\right),
\label{so6}
\end{equation}
Therefore, we consider only computation of the function $K^{(D)}\left( 
\mathbf{y,x}\right) $.

\subsection{Integral representations}

Equations (\ref{so5}) already provide an integral representation of the
solution. Since the domain of integration $S_{\infty }$ is simple, the
integrals can be written in a more explicit form. Let us consider the
Cartesian, cylindrical, and spherical coordinates of the points, 
\begin{equation}
\mathbf{r}=\left( x,y,z\right) =\left( \rho \cos \varphi ,\rho \sin \varphi
,z\right) =r\left( \sin \theta \cos \varphi ,\sin \theta \sin \varphi ,\cos
\theta \right) ,  \label{ir1}
\end{equation}%
and mark the coordinates of $\mathbf{x,y,}$ and $\mathbf{x}^{\prime }$ with
subscripts $x$,$y$, and the prime, respectively. Note then that for points $%
\mathbf{x}^{\prime }$ on the surface $z^{\prime }=0$, so 
\begin{eqnarray}
\left. \frac{\partial }{\partial z^{\prime }}\frac{1}{\left| \mathbf{x}%
^{\prime }-\mathbf{y}\right| }\right| _{z^{\prime }=0} &=&\left. \frac{z_{y}%
}{\left| \mathbf{x}^{\prime }-\mathbf{y}\right| ^{3}}\right| _{z^{\prime
}=0},  \label{ir2} \\
\left. \left| \mathbf{x}-\mathbf{x}^{\prime }\right| \right| _{z^{\prime
}=0} &=&\left. \left( \left| \mathbf{x}^{\prime }\right| ^{2}-2\mathbf{x}%
^{\prime }\cdot \mathbf{x}+\left| \mathbf{x}\right| ^{2}\right)
^{1/2}\right| _{z^{\prime }=0}=\left( \rho ^{\prime 2}-2\rho _{x}\rho
^{\prime }\cos \left( \varphi ^{\prime }-\varphi _{x}\right)
+r_{x}^{2}\right) ^{1/2}.  \notag
\end{eqnarray}%
Using the expression for the free-space Green's function (\ref{sp2}), we
obtain from Eq. (\ref{so5}) 
\begin{eqnarray}
K^{(D)}\left( \mathbf{y},\mathbf{x}\right) &=&K^{(D)}\left( \rho
_{y},\varphi _{y},z_{y};\rho _{x},\varphi _{x},z_{x}\right) =-\frac{z_{y}}{%
8\pi ^{2}}\int_{0}^{2\pi }\int_{R_{0}}^{\infty }fd\rho ^{\prime }d\varphi
^{\prime },  \label{ir3} \\
f &=&\frac{\rho ^{\prime }}{\left( \rho ^{\prime 2}-2\rho _{y}\rho ^{\prime
}\cos \left( \varphi ^{\prime }-\varphi _{y}\right) +r_{y}^{2}\right)
^{3/2}\left( \rho ^{\prime 2}-2\rho _{x}\rho ^{\prime }\cos \left( \varphi
^{\prime }-\varphi _{x}\right) +r_{x}^{2}\right) ^{1/2}}.  \notag
\end{eqnarray}%
The integrand has the following asymptotic behavior for large $\rho ^{\prime
}$ 
\begin{equation}
f=\frac{1}{\rho ^{\prime 3}}+\frac{1}{\rho ^{\prime 4}}\left( \rho _{x}\cos
\left( \varphi ^{\prime }-\varphi _{x}\right) +3\rho _{y}\cos \left( \varphi
^{\prime }-\varphi _{y}\right) \right) +O\left( \frac{r_{x}^{2}+r_{y}^{2}}{%
\rho ^{\prime 5}}\right) .  \label{ir4}
\end{equation}%
This shows that the integral converges. The infinite domain can be truncated
and we have 
\begin{equation}
K^{(D)}=K_{t}^{(D)}-\frac{z_{y}}{8\pi R_{\infty }^{2}}\left( 1+O\left( \frac{%
r_{x}^{2}+r_{y}^{2}}{R_{\infty }^{2}}\right) \right) ,\quad K_{t}^{(D)}=-%
\frac{z_{y}}{8\pi ^{2}}\int_{0}^{2\pi }\int_{R_{0}}^{R_{\infty }}fd\rho
^{\prime }d\varphi ^{\prime },  \label{ir5}
\end{equation}%
where $R_{\infty }$ is some cutoff radius. The second term of the asymptotic
expansion in Eq. (\ref{ir4}) does not contribute to the integral, since
being integrated over the angle it produces zero. A closed form for
numerical integration can be obtain by substituting the integration variable 
$\rho ^{\prime }=R_{0}/\eta $. We have then from Eq. (\ref{ir3}) 
\begin{gather}
K^{(D)}\left( \mathbf{y},\mathbf{x}\right) =-\frac{R_{0}^{2}z_{y}}{8\pi ^{2}}%
\int_{0}^{2\pi }\int_{0}^{1}\frac{\eta d\eta d\varphi ^{\prime }}{Q\left(
\eta ,\varphi ^{\prime }-\varphi _{y},\rho _{y},r_{y},R_{0}\right)
^{3}Q\left( \eta ,\varphi ^{\prime }-\varphi _{x},\rho
_{x},r_{x},R_{0}\right) },  \label{ir7} \\
Q\left( \eta ,\varphi ,\rho ,r,R\right) =\left( r^{2}\eta ^{2}-2R\rho \eta
\cos \varphi +R^{2}\right) ^{1/2}.  \notag
\end{gather}

\subsection{Spherical harmonic expansions}

While the integral representations of $K^{(D)}\left( \mathbf{y},\mathbf{x}%
\right) $ enable computation of this function for arbitrary locations of the
source and evaluation points, numerical integration can be relatively
expensive and more efficient way to compute the integral can be of practical
interest. Equation (\ref{sp6}) shows that to solve the boundary integral
equations $K^{(D)}\left( \mathbf{y},\mathbf{x}\right) $ should be evaluated
only for the case $\mathbf{x,}$ $\mathbf{y}\in \Omega _{0}$. In this case
the integral can be computed using spherical harmonic expansion. For a
single source located at $\mathbf{x}\in \Omega _{0}$ and $\mathbf{x}^{\prime
}\in \Omega _{\infty }$ the following expansion of the free space Green's
function holds \cite{Morse1953book}, \cite{Gumerov2005report}: 
\begin{equation}
\!G\left( \mathbf{x,x}^{\prime }\right) \!=\sum_{n=0}^{\infty
}\sum_{m=-n}^{n}\frac{1}{2n+1}R_{n}^{-m}\left( \mathbf{x}\right)
S_{n}^{m}\left( \mathbf{x}^{\prime }\right) ,\quad \!\left| \mathbf{x}%
^{\prime }\right| \!>R_{0}>\left| \mathbf{x}\right| .\!  \label{sh1}
\end{equation}%
Here $R_{n}^{m}\left( \mathbf{r}\right) $ and $S_{n}^{m}\left( \mathbf{r}%
\right) $ are the regular and singular (at $r=0)$ spherical basis functions 
\begin{equation}
R_{n}^{m}\left( \mathbf{r}\right) =r^{n}Y_{n}^{m}\left( \theta ,\varphi
\right) ,\quad S_{n}^{m}\left( \mathbf{r}\right) =r^{-n-1}Y_{n}^{m}\left(
\theta ,\varphi \right) ,  \label{sh2}
\end{equation}%
where $Y_{n}^{m}\left( \theta ,\varphi \right) $ are the orthonormal
spherical harmonics, 
\begin{eqnarray}
Y_{n}^{m}\left( \theta ,\varphi \right) &=&N_{n}^{m}P_{n}^{\left| m\right|
}(\xi )e^{im\varphi },\quad \xi =\cos \theta ,\quad  \label{sh3} \\
N_{n}^{m} &=&(-1)^{m}\sqrt{\frac{2n+1}{4\pi }\frac{(n-\left| m\right| )!}{%
(n+\left| m\right| )!}},\quad n=0,1,2,...,\quad m=-n,...,n.  \notag
\end{eqnarray}%
where $P_{n}^{m}\left( \xi \right) $ are the associated Legendre functions
defined by the Rodrigue's formula (see \cite{Abramowitz1964book}), 
\begin{equation}
P_{n}^{m}\left( \xi \right) =\frac{\left( -1\right) ^{m}\left( 1-\xi
^{2}\right) ^{m/2}}{2^{n}n!}\frac{d^{n+m}}{d\xi ^{n+m}}\left( 1-\xi
^{2}\right) ^{n},\quad n\geqslant 0,\quad m\geqslant 0.  \label{sh4}
\end{equation}

The derivatives of the basis functions can be expressed via the basis
functions (e.g., \cite{Gumerov2005report}) 
\begin{equation}
\frac{\partial }{\partial z}S_{n}^{m}\left( \mathbf{r}\right) =-\left(
2n+1\right) a_{n}^{m}S_{n+1}^{m}\left( \mathbf{r}\right) ,\quad  \label{sh5}
\end{equation}%
where 
\begin{eqnarray}
a_{n}^{m} &=&a_{n}^{-m}=\sqrt{\frac{(n+1+m)(n+1-m)}{\left( 2n+1\right)
\left( 2n+3\right) }},\quad \text{for }n\geqslant \left| m\right| .
\label{sh6} \\
a_{n}^{m} &=&0,\quad \text{for }n<\left| m\right| .  \notag
\end{eqnarray}%
Thus, we have 
\begin{equation}
\!\frac{\partial G\left( \mathbf{y,x}^{\prime }\right) }{\partial z^{\prime }%
}=-\sum_{n=0}^{\infty }\sum_{m=-n}^{n}a_{n}^{m}R_{n}^{-m}\left( \mathbf{y}%
\right) S_{n+1}^{m}\left( \mathbf{x}^{\prime }\right) ,\quad \!\left| 
\mathbf{x}^{\prime }\right| \!>R_{0}>\left| \mathbf{y}\right| ,\quad
z^{\prime }=z(\mathbf{x}^{\prime }).  \label{sh7}
\end{equation}%
For the computation of the integrals in (\ref{so5}) only the values for $%
z^{\prime }=0$, or $\theta ^{\prime }=\pi /2$ are needed. Using the value
for $P_{n}^{|m|}(0)$ \cite{Abramowitz1964book}, we have from (\ref{sh3}) 
\begin{eqnarray}
Y_{n}^{m}\left( \frac{\pi }{2},\varphi \right) &=&L_{n}^{m}e^{im\varphi },
\label{sh8} \\
L_{n}^{m} &=&N_{n}^{m}P_{n}^{|m|}(0)=l_{n}^{m}\frac{\left( -1\right)
^{\left( n+\left| m\right| \right) /2}\left( n+\left| m\right| \right) !}{%
2^{n}\left( \frac{n-\left| m\right| }{2}\right) !\left( \frac{n+\left|
m\right| }{2}\right) !}N_{n}^{m},  \notag \\
l_{n}^{m} &=&\left\{ 
\begin{array}{c}
0,\quad n+m=2l+1,\quad l=0,1,... \\ 
1,\quad n+m=2l,\quad l=0,1,...%
\end{array}%
\right. .  \notag
\end{eqnarray}%
Hence, 
\begin{equation}
\left. S_{n}^{m}\left( \mathbf{r}\right) \right| _{\mathbf{r\in }S_{\infty
}}=L_{n}^{m}\rho ^{-n-1}e^{im\varphi },\quad \left( \left. \rho \right| _{%
\mathbf{r\in }S_{\infty }}=\left. r\sin \theta \right| _{\mathbf{r\in }%
S_{\infty }}=\left. r\right| _{\mathbf{r\in }S_{\infty }}\right)  \label{sh9}
\end{equation}%
This shows that 
\begin{eqnarray}
\int_{S_{\infty }}S_{n^{\prime }+1}^{m^{\prime }}\left( \mathbf{r}\right)
S_{n}^{m}\left( \mathbf{r}\right) dS\left( \mathbf{r}\right) &=&L_{n^{\prime
}+1}^{m^{\prime }}L_{n}^{m}\int_{R}^{\infty }\rho ^{-n-n^{\prime }-3}\rho
d\rho \int_{0}^{2\pi }e^{i\left( m+m^{\prime }\right) \varphi }d\varphi
\label{sh10} \\
&=&2\pi \delta _{m^{\prime },-m}\frac{L_{n^{\prime }+1}^{m^{\prime
}}L_{n}^{m}}{n+n^{\prime }+1}R^{-n-n^{\prime }-1}  \notag
\end{eqnarray}%
where $\delta _{m^{\prime }m}$ is the Kronecker delta.

Using Eqs (\ref{sh1}), (\ref{sh7}), (\ref{sh9}), and (\ref{sh10}) we obtain
from Eq. (\ref{so5}) 
\begin{eqnarray}
K^{(D)}\left( \mathbf{y,x}\right) &=&\sum_{n=0}^{\infty
}\sum_{m=-n}^{n}\sum_{n^{\prime }=0}^{\infty }\sum_{m^{\prime }=-n^{\prime
}}^{n^{\prime }}\frac{2a_{n^{\prime }}^{m^{\prime }}}{2n+1}R_{n}^{-m}\left( 
\mathbf{x}\right) R_{n^{\prime }}^{-m^{\prime }}\left( \mathbf{y}\right)
\int_{S_{\infty }}S_{n^{\prime }+1}^{m^{\prime }}\left( \mathbf{x}^{\prime
}\right) S_{n}^{m}\left( \mathbf{x}^{\prime }\right) dS\left( \mathbf{x}%
^{\prime }\right)  \label{sh11} \\
&=&\sum_{n=0}^{\infty }\sum_{m=-n}^{n}\sum_{n^{\prime }=0}^{\infty
}I_{n^{\prime }n}^{m}R_{n}^{-m}\left( \mathbf{x}\right) R_{n^{\prime
}}^{m}\left( \mathbf{y}\right) =\sum_{m=-\infty }^{\infty }\sum_{n=\left|
m\right| +1}^{\infty }\sum_{n^{\prime }=\left| m\right| }^{\infty
}I_{nn^{\prime }}^{m}R_{n^{\prime }}^{-m}\left( \mathbf{x}\right)
R_{n}^{m}\left( \mathbf{y}\right)  \notag
\end{eqnarray}%
where 
\begin{equation}
I_{nn^{\prime }}^{m}=J_{nn^{\prime }}^{m}R^{-n-n^{\prime }-1},\quad
J_{nn^{\prime }}^{m}=\frac{4\pi a_{n}^{m}L_{n+1}^{m}L_{n^{\prime }}^{m}}{%
\left( 2n^{\prime }+1\right) \left( n+n^{\prime }+1\right) }.  \label{sh12}
\end{equation}%
Here we used symmetries $a_{n}^{m}=a_{n}^{-m}$ and $L_{n}^{m}=L_{n}^{-m}$.
Summation in the last sum of Eq (\ref{sh11}) can start at $n=\left| m\right| 
$, but at this value $L_{n+1}^{m}=0$ (see Eq. (\ref{sh8})). Note that $%
I_{nn^{\prime }}^{m}=I_{nn^{\prime }}^{-m}$, but the coefficients are
non-symmetric with respect to the lower indices, $I_{nn^{\prime }}^{m}\neq
I_{n^{\prime }n}^{m}$. Indeed, we have 
\begin{equation}
I_{n^{\prime }n}^{m}I_{nn^{\prime }}^{m}=\frac{\left( 4\pi \right)
^{2}a_{n}^{m}a_{n^{\prime }}^{m}R^{-2n-2n^{\prime }-2}}{\left( 2n+1\right)
\left( 2n^{\prime }+1\right) \left( n+n^{\prime }+1\right) ^{2}}%
L_{n}^{m}L_{n+1}^{m}L_{n^{\prime }}^{m}L_{n^{\prime }+1}^{m}=0.
\label{sh12.1}
\end{equation}%
This is due to the fact $L_{n}^{m}L_{n+1}^{m}=0$ (see definition (\ref{sh8}%
)). Hence, if at some $n^{\prime }$ and $n$ we have $I_{n^{\prime
}n}^{m}\neq 0$ (and this is true as not all $I_{n^{\prime }n}^{m}$ are
zeros), then due to Eq. (\ref{sh12.1}) $I_{nn^{\prime }}^{m}=0$, means $%
I_{n^{\prime }n}^{m}\neq $ $I_{nn^{\prime }}^{m}$. Non-symmetry of the
coefficients also shows that $K^{(D)}\left( \mathbf{y,x}\right) $ is
non-symmetric, $K^{(D)}\left( \mathbf{y,x}\right) \neq K^{(D)}\left( \mathbf{%
x,y}\right) $ since the spherical harmonics form a complete orthogonal basis
on a sphere (and $R_{n}^{m}\left( \mathbf{y}\right) $ are proportional to
them).

\section{Numerical methods}

\subsection{Dimensionless forms}

Since the Laplace equation does not have intrinsic length scales, we have
for arbitrary $R$ (see Eqs (\ref{ir7}) and (\ref{sh11})) 
\begin{eqnarray}
K^{(D)}\left( \mathbf{y},\mathbf{x;}R\right) &=&\frac{1}{R}\widetilde{K}%
^{(D)}\left( \frac{\mathbf{y}}{R},\frac{\mathbf{x}}{R}\right) ,  \label{nm1}
\\
\widetilde{K}^{(D)}\left( \mathbf{y},\mathbf{x}\right) &=&-\frac{z_{y}}{8\pi
^{2}}\int_{0}^{2\pi }\int_{0}^{1}\frac{\eta d\eta d\varphi ^{\prime }}{%
Q\left( \eta ,\varphi ^{\prime }-\varphi _{y},\rho _{y},r_{y},1\right)
^{3}Q\left( \eta ,\varphi ^{\prime }-\varphi _{x},\rho _{x},r_{x},1\right) },
\notag \\
\widetilde{K}^{(D)}\left( \mathbf{y},\mathbf{x}\right) &=&\sum_{n=0}^{\infty
}\sum_{m=-n}^{n}U_{n}^{m}\left( \mathbf{x}\right) R_{n}^{m}\left( \mathbf{y}%
\right) ,\quad U_{n}^{m}\left( \mathbf{x}\right) =\sum_{n^{\prime }=\left|
m\right| }^{\infty }J_{nn^{\prime }}^{m}R_{n^{\prime }}^{-m}\left( \mathbf{x}%
\right) ,\quad \left( \left| \mathbf{y}\right| <1,\left| \mathbf{x}\right|
<1\right) .  \notag
\end{eqnarray}%
where $\widetilde{K}^{(D)}$ is the dimensionless function, which does not
depend on parameter $R$, which is seen from its integral representation. We
put an extra argument $R$ to the dimensional function $K^{(D)}$ to show its
parametric dependence on $R$. The series for $\widetilde{K}^{(D)}\left( 
\mathbf{y},\mathbf{x}\right) $ and $U_{n}^{m}\left( \mathbf{x}\right) $%
converge for $\mathbf{x}$ and $\mathbf{y}$ located inside the unit ball.

\subsection{Improvement of convergence}

Our numerical experiments show that computation of the kernel function via
double integral representation (\ref{nm1}) provides accurate results, but
such computations are somehow slow as they require evaluation of integrals
for every pair $\mathbf{x}$ and $\mathbf{y}$. On the other hand, recursive
computations of the spherical basis functions provide a fast procedure.

For practical use all infinite series should be truncated with some
truncation number $p$ ($n<p,$ $n^{\prime }<p$), which should be selected
based on the error bounds. The radius of convergence of series (\ref{nm1})
is $1$ both, with respect to $\mathbf{x}$ and $\mathbf{y}$. This means that
for points in $\Omega _{0}$ located closely to the boundary with domain $%
\Omega _{\infty }$ the truncation number can be large, which may reduce the
efficiency of computations.

Several tricks can be proposed to improve the convergence. First, one can
simply implement integral and series representations of the kernel and
switch between algorithms based on the values $\left| \mathbf{y}\right| $
and $\left| \mathbf{x}\right| $. Second, we can effectively reduce the size
of the domain, where $K^{(D)}\left( \mathbf{y},\mathbf{x}\right) $ should be
computed. This is based on the observation that for the evaluation points
located on plane $z=0$ we have from Eqs (\ref{sh2}), (\ref{sh3}), (\ref{sh8}%
), and (\ref{sh12}) 
\begin{equation}
J_{nn^{\prime }}^{m}\left. R_{n}^{m}\left( \mathbf{y}\right) \right|
_{z_{y}=0}=\frac{4\pi a_{n}^{m}L_{n^{\prime }}^{m}}{\left( 2n^{\prime
}+1\right) \left( n+n^{\prime }+1\right) }\left( L_{n+1}^{m}L_{n}^{m}\right)
\rho _{y}^{n}e^{-im\varphi _{y}}=0,  \label{nm5}
\end{equation}%
due to $L_{n+1}^{m}L_{n}^{m}=0$. Hence, 
\begin{equation}
K^{(D)}\left( \mathbf{y},\mathbf{x;}R\right) =0,\quad z_{y}=0.  \label{nm6}
\end{equation}%
Also, relation (\ref{nm6}) follows from the integral representation (\ref%
{ir3}).

Consider now an extended computational domain, which can be denoted as $%
\Omega _{e}$, and which includes $\Omega _{0}$ and the part of $%
\Omega_{\infty }$ located between the concentric spheres of radii $R_{0}$
and $R_{e}>R_{0}.$ Respectively, $S_{ge}$ denotes the part of ground surface 
$S_{g}$ in domain $\Omega _{e}$ (see Fig. \ref{Fig1}). All formulae derived
hold in case if we replace $\Omega _{0}$ with $\Omega _{e}$, $R_{0}$ with $%
R_{e}$ and $S_{g0}$ with $S_{ge}$. So, instead of Eq. (\ref{sp5}) we should
have 
\begin{equation}
\phi \left( \mathbf{y}\right) =\int_{S\cup S_{ge}}\sigma _{e}^{(D)}\left(%
\mathbf{x}\right) \left[ G\left( \mathbf{y},\mathbf{x}\right) +K^{(D)}\left(%
\mathbf{y},\mathbf{x;}R_{e}\right) \right] dS\left( \mathbf{x}\right).
\label{nm7}
\end{equation}
The charge density $\sigma _{e}^{(D)}\left( \mathbf{x}\right) $ can be
determined from the boundary conditions 
\begin{eqnarray}
\int_{S\cup S_{ge}}\sigma _{e}^{(D)}\left( \mathbf{x}\right) \left[ G\left(%
\mathbf{y},\mathbf{x}\right) +K^{(D)}\left( \mathbf{y},\mathbf{x;}%
R_{e}\right) \right] dS\left( \mathbf{x}\right) &=&\left\{ 
\begin{array}{c}
V\left( \mathbf{y}\right) ,\quad \mathbf{y}\in S, \\ 
0,\quad \mathbf{y}\in S_{g0},%
\end{array}%
\right. ,  \label{nm8} \\
\int_{S\cup S_{ge}}\sigma _{e}^{(D)}\left( \mathbf{x}\right) G\left( \mathbf{%
y},\mathbf{x}\right) dS\left( \mathbf{x}\right) &=&0,\quad \mathbf{y}\in
S_{ge}\backslash S_{g0}.  \notag
\end{eqnarray}
The latter holds since $K^{(D)}\left( \mathbf{y},\mathbf{x;}R_{e}\right) =0,$
for $\mathbf{y}\in S_{ge}\backslash S_{g0}$. Hence, both for the field
computations (\ref{nm7}) and determination of the charge density only values
of $K^{(D)}\left( \mathbf{y},\mathbf{x;}R_{e}\right) $ for $\mathbf{y}\in
\Omega _{0}$\textbf{\ }are needed. On the other hand, $K^{(D)}\left( \mathbf{%
y},\mathbf{x;}R_{e}\right) $ should be computed for $\mathbf{x}\in
\Omega_{e} $ with a notice that for $\mathbf{x}\in S_{ge}\backslash S_{g0}$
we have $z_{x}=0.$ We developed an efficient recursive procedure for such
computations described in the next subsection.

\subsection{Computations for the ground points on the extended domain}

So, the problem for use of the extended domain is to compute function $%
U_{n}^{m}\left( \mathbf{x}\right) $ represented in Eq. (\ref{nm1}) by
infinite series at $z_{x}$ $=0$, at $\xi \in \left[ R_{0}/R_{e},1\right) ,$%
where $\xi =\left| \mathbf{x}\right| $. Integral representation of $%
U_{n}^{m}\left( \mathbf{x}\right) $ can be obtained by substituting of Eq. (%
\ref{sh7}) into integral (\ref{so5}). In the dimensionless form we have 
\begin{equation}
\widetilde{K}^{(D)}\left( \mathbf{y,x}\right) =\sum_{n=0}^{\infty
}\sum_{m=-n}^{n}U_{n}^{m}\left( \mathbf{x}\right) R_{n}^{m}\left( \mathbf{y}%
\right) =\sum_{n=0}^{\infty }\sum_{m=-n}^{n}2a_{n}^{m}R_{n}^{m}\left( 
\mathbf{y}\right) \int_{1}^{\infty }\int_{0}^{2\pi }G\left( \mathbf{x},%
\mathbf{x}^{\prime }\right) S_{n+1}^{-m}\left( \mathbf{x}^{\prime }\right)
dS\left( \mathbf{x}^{\prime }\right) .\quad  \label{ed1}
\end{equation}%
So, 
\begin{equation}
\left. U_{n}^{m}\left( \mathbf{x}\right) \right| _{z_{x}=0}=\frac{1}{2\pi }%
a_{n}^{m}L_{n+1}^{m}e^{-im\varphi _{x}}u_{n}^{m}\left( \xi \right) ,\quad
\xi =\rho _{x},  \label{ed2}
\end{equation}%
\begin{eqnarray}
u_{n}^{m}\left( \xi \right) &=&\int_{1}^{\infty }\frac{d\eta }{\eta ^{n+1}}%
\int_{0}^{2\pi }\frac{e^{im\varphi }d\varphi }{\left( \eta ^{2}-2\xi \eta
\cos \varphi +\xi ^{2}\right) ^{1/2}}=\frac{1}{\xi ^{n+1}}\int_{0}^{\xi
}\zeta ^{n}w_{m}\left( \zeta \right) d\zeta ,  \label{ed3} \\
w_{m}\left( \xi \right) &=&\int_{0}^{2\pi }\frac{e^{im\varphi }d\varphi }{%
\left( 1-2\xi \cos \varphi +\xi ^{2}\right) ^{1/2}},\quad u_{n}^{-m}\left(
\xi \right) =u_{n}^{m}\left( \xi \right) ,\quad w_{-m}\left( \xi \right)
=w_{m}\left( \xi \right) .  \notag
\end{eqnarray}%
These functions can be expressed via the complete elliptic integrals of the
first and second kinds, $K$ and $E$ (see\cite{Abramowitz1964book}), and
recursive relations (to shorten notation we drop the argument $\xi $ for
functions $u_{n}^{m}\left( \xi \right) $ and $w_{m}\left( \xi \right) $) 
\begin{eqnarray}
u_{n}^{0} &=&\frac{1}{n^{2}\xi ^{2}}\left[ 4E\left( \xi ^{2}\right)
-4n\left( 1-\xi ^{2}\right) K\left( \xi ^{2}\right) +\left( n-1\right)
^{2}u_{n-2}^{0}\right] ,\quad n=1,3,5,...  \label{ed4.1} \\
u_{n+1}^{1} &=&\frac{1}{\left( 2n+3\right) \xi }\left[ \left( n+1\right)
u_{n}^{0}+\left( n+2\right) \xi ^{2}u_{n+2}^{0}+4\left( 1-\xi ^{2}\right)
K\left( \xi ^{2}\right) -8E\left( \xi ^{2}\right) \right] ,\quad n=1,3,5,...
\label{ed4.2} \\
u_{n+1}^{m+1} &=&\frac{2}{\left( 2n+3\right) \xi }\left[ \left( n+1\right)
u_{n}^{m}+\left( n+2\right) \xi ^{2}u_{n+2}^{m}-v_{m}\right] -u_{n+1}^{m-1},%
\text{ }m=1,2,...,\,n=m+1,m+3,...  \label{ed4.3} \\
v_{m} &=&\left( 1+\xi ^{2}\right) w_{m}-\xi \left( w_{m+1}+w_{m-1}\right)
,\quad m=1,2,...  \label{ed4.4} \\
w_{m} &=&\frac{1}{\xi }\left( 1+\xi ^{2}\right) \frac{2m-2}{2m-1}w_{m-1}-%
\frac{2m-3}{2m-1}w_{m-2},\quad m=2,3,...,  \label{ed4.5} \\
w_{0} &=&4K\left( \xi ^{2}\right) ,\quad w_{1}=\frac{4}{\xi }\left[ K\left(
\xi ^{2}\right) -E\left( \xi ^{2}\right) \right] ,  \label{ed4.6} \\
K\left( \mu \right) &=&\int_{0}^{\pi /2}\left( 1-\mu \sin ^{2}\varphi
\right) ^{-1/2}d\varphi ,\quad E\left( \mu \right) =\int_{0}^{\pi /2}\left(
1-\mu \sin ^{2}\varphi \right) ^{1/2}d\varphi .  \label{ed4.7}
\end{eqnarray}%
Appendix A describes how these relations are obtained. We remark that
because $L_{n+1}^{m}=0$ when $n+m$ is an even number, to obtain $U_{n}^{m}$
in Eq. (\ref{ed2}) we need to compute functions $u_{n}^{m}$ only at odd
values of $n+m$, which is reflected in the above recurrences. The recurrence
for $u_{n}^{m}$ works in layers in subscript $n$. So, first, we obtain $%
u_{n}^{0}$ for all $n$ needed, then $u_{n}^{1}$, and so on until $u_{n}^{m}$
at maximum $m$ needed. For truncation number $p$, we have max$\left(
n\right) =p-1$, max$\left( m\right) =p-2$. As recurrence (\ref{ed4.3}) show
that the maximum computed $n$ at $m=p-3$ should be $p$, while at $m=0$ it
should be $2p-3$. It is also remarkable that recurrence (\ref{ed4.1}) can be
used first at $n=1$, which provides the value of $u_{1}^{0}$, and then
applied for other $u_{n}^{0}$ in the sequence.

\subsection{Boundary elements}

Equation (\ref{nm8}) can be solved using the method of moments (MoM) or
collocation boundary element methods (BEM). In the present study we use the
BEM. The surface is discretized by triangular elements, the integrals
involving Green's function are computed analytically assuming that the
charge density is constant over the elements (constant panel approximation).
The integrals involving kernel $K^{(D)}\left( \mathbf{y},\mathbf{x;}%
R_{c}\right) $ were computed using center panel quadrature, 
\begin{equation}
\int_{\Delta _{j}}\sigma \left( \mathbf{x}\right) K^{(D)}\left( \mathbf{y},%
\mathbf{x;}R_{c}\right) dS\left( \mathbf{x}\right) \approx \sigma
_{j}w_{j}K^{(D)}\left( \mathbf{y},\mathbf{x}_{j}^{c}\mathbf{;}R_{c}\right) ,
\label{nm9}
\end{equation}%
where $\sigma _{j}$ is the mean value of the charge density, while $w_{j}$
and $\mathbf{x}_{j}^{c}$ are the area and the centroid of triangle $\Delta
_{j}$. In the standard BEM realization the integrals over the kernels with
Green's functions were computed analytically. 
\begin{equation}
\int_{\Delta _{j}}\sigma \left( \mathbf{x}\right) G\left( \mathbf{y},\mathbf{%
x}\right) dS\left( \mathbf{x}\right) \approx \sigma _{j}L_{j}\left( \mathbf{y%
}\right) ,  \label{nm10}
\end{equation}%
\begin{eqnarray}
L_{j}\left( \mathbf{y}\right) &=&\int_{\Delta _{j}}G\left( \mathbf{y},%
\mathbf{x}\right) dS\left( \mathbf{x}\right) =\sum_{q=1}^{3}\left[ H\left(
l_{jq}-x_{jq},h_{j},z_{jq}\right) -H\left( -x_{jq},h_{j},z_{jq}\right) %
\right] ,  \label{nm11} \\
x_{jq} &=&\left( \mathbf{y-x}_{jq}\right) \cdot \mathbf{i}_{jq},\quad
h_{j}=\left| \left( \mathbf{y-x}_{j1}\right) \cdot \mathbf{n}_{j}\right|
,\quad z_{jq}=\left( \mathbf{y-x}_{jq}\right) \cdot \mathbf{n}_{jq},  \notag
\\
\mathbf{i}_{jq} &=&\frac{1}{l_{jq}}\left( \mathbf{x}_{j,q(\func{mod}3)+1}-%
\mathbf{x}_{jq}\right) ,\quad l_{jq}=\left| \mathbf{x}_{j,q(\func{mod}3)+1}-%
\mathbf{x}_{jq}\right| ,\quad \mathbf{n}_{jq}=\mathbf{i}_{jq}\times \mathbf{n%
}_{j},  \notag \\
\quad H\left( x,y,z\right) &=&y\left( \arctan \frac{x}{z}-\arctan \frac{yx}{%
zr}\right) -z\ln \left| r+x\right| ,\quad r=\sqrt{x^{2}+y^{2}+z^{2}},  \notag
\end{eqnarray}%
where $\mathbf{n}_{j}$ is the normal to triangle $\Delta _{j}$ directed into
the computational domain.

\subsection{Fast multipole method}

For large problems the BEM can be accelerated using the fast multipole
method (FMM) (e.g., \cite{Adelman2017IEEE}). As shown, the kernel $%
K^{(D)}\left( \mathbf{y},\mathbf{x;}R_{c}\right) $ is either zero or can be
factored for any pair of points inside the computational domain of radius $%
R_{c}$. This shows that boundary integral equation (\ref{nm8}) can be
written as 
\begin{gather}
\int_{S\cup S_{ge}}\!\!\!\!\!\!\!\!\!\!\sigma _{e}^{(D)}\left( \mathbf{x}%
\right) G\left( \mathbf{y},\mathbf{x}\right) dS\left( \mathbf{x}\right) +
\label{nm12} \\
\frac{1}{R_{e}}\sum_{n=0}^{p-1}\!\sum_{m=-n}^{n}R_{n}^{m}\left( \frac{%
\mathbf{y}}{R_{e}}\right) \int_{S\cup S_{ge}}\!\!\!\!\!\!\!\!\!\!\sigma
_{e}^{(D)}\left( \mathbf{x}\right) U_{n}^{m}\left( \frac{\mathbf{x}}{R_{e}}%
\right) dS\left( \mathbf{x}\right) =\left\{ 
\begin{array}{c}
V\left( \mathbf{y}\right) ,\quad \mathbf{y}\in S, \\ 
0,\quad \mathbf{y}\in S_{ge},%
\end{array}%
\right.  \notag
\end{gather}%
where the infinite sum over $n$ is truncated. The first integral in the left
hand side can be computed using the FMM with the free space Green's
function. In our implementation, we used some nearfield radius $r_{nf}$. So,
for pairs of points $\left( \mathbf{y},\mathbf{x}\right) $ located closer
than $r_{nf}$ the integrals over triangles were computed analytically using
Eqs (\ref{nm10}) and (\ref{nm11}), while for larger distances the integrals
were computed using center point approximation, similarly to (\ref{nm9}).

The most important feature of the global kernel factorization is that the
number of operations to compute the matrix-vector product for the kernel is $%
O\left( p^{2}N\right) $, not $O\left( N^{2}\right) $, since all integrals
over $\mathbf{x}$ can be computed independently from the evaluation points $%
\mathbf{y}$ ($N$ is the number of boundary elements). So, computation of
this sum has the asymptotic complexity consistent with the FMM, and the
method can be applied for solution of large problems, where the overall
linear system is solved iteratively.

\subsection{Real basis functions}

Kernel $\widetilde{K}^{(D)}$ is real, but is presented in Eq. (\ref{nm1})
via sums of complex valued functions $R_{n}^{m}$. A fast procedure can be
developed based on recursive computations of real basis functions (see \cite%
{Gumerov2008JCP}), 
\begin{equation}
\underline{R_{n}^{m}}\left( \mathbf{r}\right) =\frac{\left( -1\right) ^{n+m}%
}{(n+\left| m\right| )!}r^{n}P_{n}^{\left| m\right| }(\cos \theta )\left\{ 
\begin{array}{c}
\cos m\varphi ,\quad m\geqslant 0, \\ 
\sin m\varphi ,\quad m<0,%
\end{array}%
\right.  \label{nm2}
\end{equation}%
which can be computed recursively as 
\begin{eqnarray}
\underline{R_{0}^{0}} &=&1,\quad \underline{R_{1}^{1}}=-\frac{1}{2}x,\quad 
\underline{R_{1}^{-1}}=\frac{1}{2}y,  \label{nm3} \\
\underline{R_{m}^{m}} &=&-\frac{1}{2m}\left( x\underline{R_{m-1}^{m-1}}+y%
\underline{R_{m-1}^{-\left( m-1\right) }}\right) ,\quad m=2,3,...,  \notag \\
\quad \underline{R_{m}^{-m}} &=&\frac{1}{2m}\left( y\underline{R_{m-1}^{m-1}}%
-x\underline{R_{m-1}^{-\left( m-1\right) }}\right) ,\quad m=2,3,...,  \notag
\\
\underline{R_{m+1}^{\pm m}} &=&-z\underline{R_{m}^{\pm m}},\quad m=0,1,..., 
\notag \\
\underline{R_{n+1}^{\pm m}} &=&-\frac{\left( 2n+1\right) z\underline{%
R_{n}^{\pm m}}+r^{2}\underline{R_{n-1}^{\pm m}}}{\left( n+1\right) ^{2}-m^{2}%
},\quad m=0,1,...,\text{ }n=m+1,m+2,....\quad  \notag
\end{eqnarray}

Equations (\ref{nm1}), (\ref{sh2}), (\ref{sh3}), (\ref{sh6}), (\ref{sh8}%
),and (\ref{sh12}) show then that 
\begin{eqnarray}
\widetilde{K}^{(D)}\left( \mathbf{y},\mathbf{x}\right)
&=&\sum_{n=0}^{\infty}\sum_{m=-n}^{n}\underline{U_{n}^{m}}\left( \mathbf{x}%
\right) \underline{R_{n}^{m}}\left( \mathbf{y}\right) ,\quad \underline{%
U_{n}^{m}}\left(\mathbf{x}\right) =-\frac{2-\delta _{m,0}}{4\pi }%
\nu_{n+1}^{m}\sum_{n^{\prime }=\left| m\right| }^{\infty }\frac{\nu
_{n^{\prime}}^{m}}{n^{\prime }+n+1}\underline{R_{n^{\prime }}^{m}}\left( 
\mathbf{x}\right) ,  \label{nm4} \\
\nu _{n}^{m} &=&l_{n}^{m}\left( -1\right) ^{\left( n+\left| m\right|
\right)/2}(n-\left| m\right| -1)!!(n+\left| m\right| -1)!!.  \notag
\end{eqnarray}
Also, using Eq. (\ref{ed2}) we have for points $z_{x}=0$ 
\begin{equation}
\underline{U_{n}^{m}}\left( \mathbf{x}\right) =-\frac{2-\delta _{m,0}}{%
8\pi^{2}}\nu _{n+1}^{m}u_{n}^{m}\left( \xi \right) \left\{%
\begin{array}{c}
\cos \left( m\varphi _{x}\right) ,\quad m\geqslant 0, \\ 
\sin \left( m\varphi _{x}\right) ,\quad m<0.%
\end{array}%
\right.  \label{nm4.1}
\end{equation}

\section{Complexity and optimization}

\subsection{Kernel computations}

Assume there are $N$ sources and $M$ receivers in domain $\Omega _{0}$.
Truncation of the series representation of $K^{(D)}\left( \mathbf{y},\mathbf{%
x;}R_{e}\right) $ for $\mathbf{y\in }$ $\Omega _{0}$ with truncation number $%
p$ produces errors $\epsilon \sim (R_{0}/R_{e})^{p}$. If there are $N_{e}$
sources in the extended domain $\Omega _{e}$ then the computational cost for
factored $K^{(D)}\left( \mathbf{y},\mathbf{x;}R_{e}\right) $ can be
estimated as 
\begin{equation}
C_{fact}=O\left( Mp^{2}\right) +O\left( Np^{3}\right) +O\left(
(N_{e}-N)p^{2}\right) ,  \label{co1}
\end{equation}%
where the first term is the cost to compute $p^{2}$ functions \underline{$%
R_{n}^{m}$}$\left( \mathbf{y/}R_{e}\right) $, the second term is the cost to
compute $p^{2}$ functions \underline{$U_{n}^{m}$}$\left( \mathbf{x/}%
R_{e}\right) $ for $\mathbf{x}\in \Omega _{0}$ using series (\ref{nm4}), and
the third term is the cost to compute $p^{2}$ functions \underline{$%
U_{n}^{m} $}$\left( \mathbf{x/}R_{e}\right) $ for $\mathbf{x}\in \Omega _{e}$
\TEXTsymbol{\backslash}$\Omega _{0}$ $\left( z_{x}=0\right) $ using Eq. (\ref%
{nm4.1}) and recursions (\ref{ed4.1})-(\ref{ed4.7}) for $u_{n}^{m}\left( \xi
\right) $. Assume now that the sources are uniformly distributed in domain $%
\Omega _{e}$\TEXTsymbol{\backslash}$\Omega _{0}$, with the density $\rho $
points per surface area. In this case $N_{e}-N=\pi \rho \left(
R_{e}^{2}-R_{0}^{2}\right) $. For prescribed accuracy $\epsilon $ we have $%
R_{e}=R_{0}\epsilon ^{-1/p}$, which brings the following expression for the
computational cost 
\begin{equation}
C_{fact}=AMp^{2}+BNp^{3}+C\pi \rho R_{0}^{2}\left( \epsilon ^{-2/p}-1\right)
p^{2},  \label{co2}
\end{equation}%
where $A,B$, and $C$ are some asymptotic constants.

Considering $C_{fact}$ as a function of $p$, we can see that 
\begin{equation}
\frac{dC_{fact}}{dp}=2AMp+3BNp^{2}+2C\pi \rho R_{0}^{2}p\left[ \left(
\epsilon ^{-2/p}-1\right) -\epsilon ^{-2/p}\frac{1}{p}\ln \frac{1}{\epsilon }%
\right] .  \label{co3}
\end{equation}%
The term in the square brackets has a single zero at $p=p_{c}\left( \epsilon
\right) $ and is positive at 
\begin{equation}
p>p_{c}\left( \epsilon \right) =\frac{1}{\alpha }\ln \frac{1}{\epsilon }%
,\quad \alpha =0.79681213...\approx 0.8.  \label{co4}
\end{equation}%
So, function $C_{fact}\left( p\right) $ grows monotonically at $%
p>p_{c}\left( \epsilon \right) $ for any values of $M,N,$ and $\pi \rho
R_{0}^{2}.$ The function can grow or decay at $p<p_{c}\left( \epsilon
\right) $ depending on $M,N,$ and $\pi \rho R_{0}^{2}$. Note now, that $%
p=p_{c}\left( \epsilon \right) $ corresponds to 
\begin{equation}
\left( \frac{R_{e}}{R_{0}}\right) _{c}=\epsilon ^{-1/p_{c}}=e^{\alpha
}=\beta \approx 2.2.  \label{co5}
\end{equation}%
This means, that $C_{fact}$ being considered as a function of $R_{e}$ decays
in the region $R_{0}<R_{e}\leqslant \beta R_{0}$ and may continue to decay
or start to grow at $R_{e}>\beta R_{0}$. It is remarkable, that $\beta $
does not depend on $\epsilon $, so for computations using the factored
kernel the optimal radius of the extended domain should be always larger
than $\beta R_{0}$. Also, it can be noted that the values of $p_{c}\left(
\epsilon \right) $ are rather small. For example, at $\epsilon =10^{-4}$ we
have $p_{c}\left( \epsilon \right) \approx $ $11.6$, which means that $p=12$
should be sufficient to achieve the required accuracy at $R_{e}=\beta R_{0}$%
. This also shows that at large enough $N$ and $M$ computational costs in
Eq.(\ref{co2}) are much smaller than that for kernel computations using the
integral representation, which formally is $O\left( NM\right) $ with the
asymptotic constant that can be substantially large as the cost of numerical
quadrature also depends on the required accuracy. 
\begin{figure}[tbh]
\vspace{-20pt}
\par
\begin{center}
\includegraphics[width=0.9\textwidth, trim=0 0.9in 0.5in -0.25in]{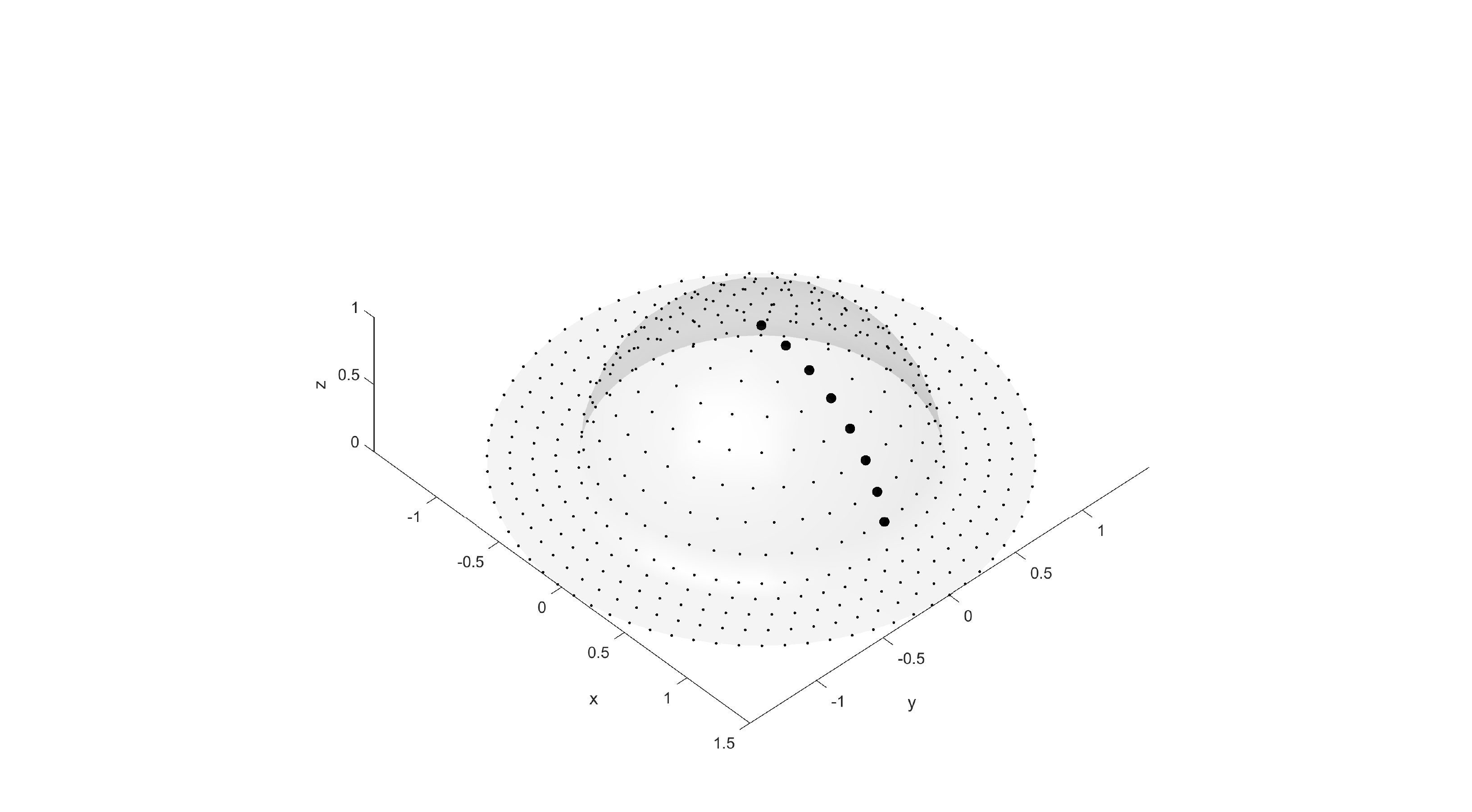}
\end{center}
\caption{The distribution of the sources (the small dots) and receivers (the
bold dots) for the kernel accuracy and computational complexity tests. }
\label{Fig2}
\end{figure}

\begin{figure}[tbh]

\par
\begin{center}
\includegraphics[width=0.9\textwidth, trim=0.5in .25in 0.5in 0.25in]{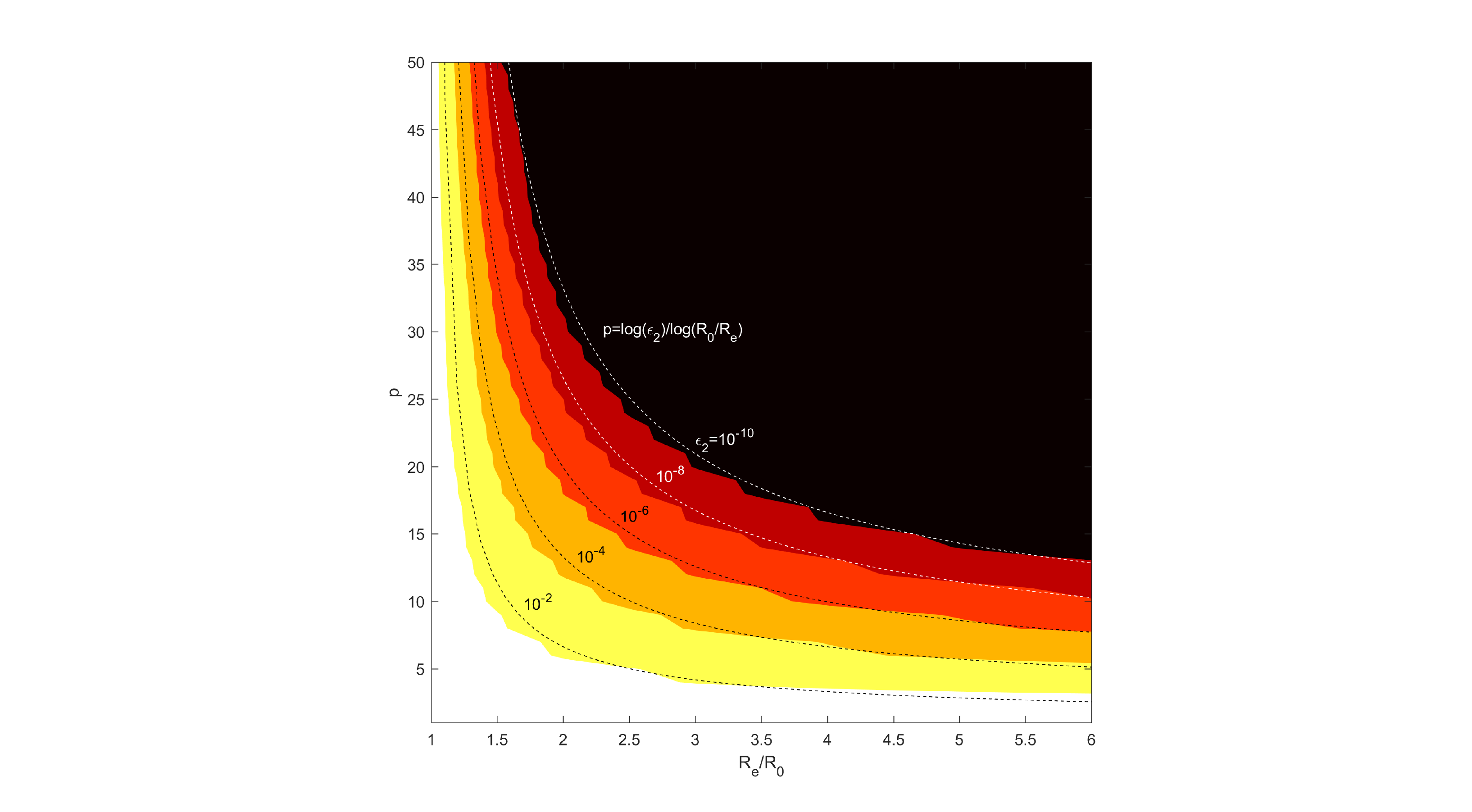}
\end{center}
\caption{The computed domains (colored) in $\left( R_{e}/R_{0},p\right) $%
-plane corresponding to different values of the relative $L_{2}$-norm error $%
\protect\epsilon _{2}$ for the factored kernel representation (Eqs (\ref{nm4}%
) and (\ref{nm4.1})) using truncation number $p$. The reference solution is
obtained using the integral representation and the Matlab integral2
procedure with the absolute and relative tolerance 10$^{-12}$. The
boundaries of the domains correlate well with dependences, $p=\log \protect%
\epsilon _{2}/\log \left( R_{0}/R_{e}\right) $ shown by the dashed curves.
The values of $\protect\epsilon _{2}$ are shown near the curves. }
\label{Fig3}
\end{figure}

\begin{figure}[tbh]
\par
\begin{center}
\includegraphics[width=0.9\textwidth, trim=0.5in .25in 0.5in 0.25in]{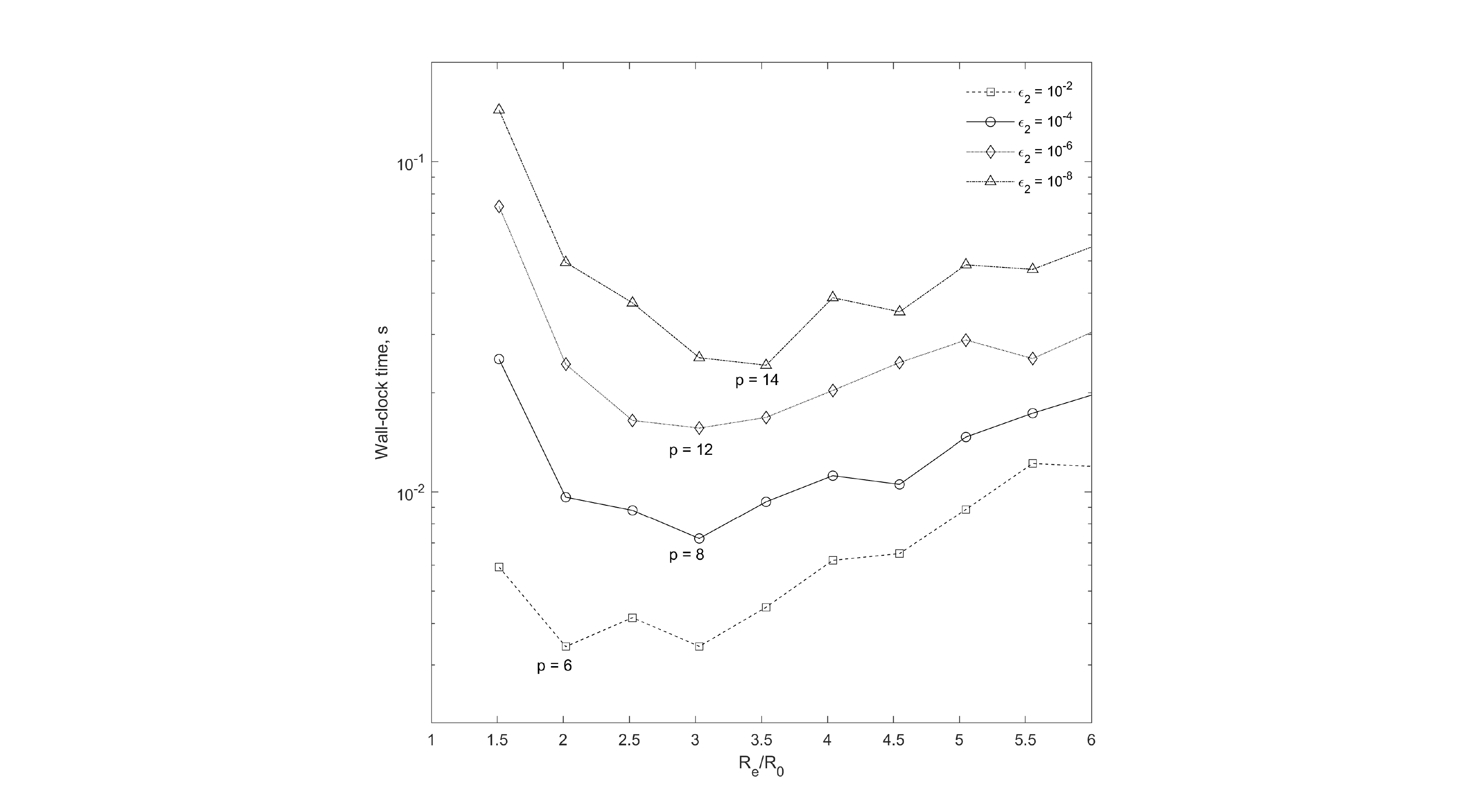}
\end{center}
\caption{The minimal wall-clock times (in seconds, Matlab implementation)
required to compute the factored kernel representation (Eqs (\ref{nm4}) and (%
\ref{nm4.1})) to achieve the presciribed error $\protect\epsilon _{2}$ as a
function of $R_{e}/R_{0}$. The values of truncation numbers for optimum
domain sizes are shown near the minima of the curves. }
\label{Fig4}
\end{figure}

To validate these findings we conducted some numerical experiments. First,
we did error estimates of the factored solutions using comparisons with the
error-controlled integral computations (Matlab integral2 function). For this
purpose due to the symmetry of the problem, all receivers can be placed on
the intersection of a unit hemisphere (``bump'') with plane $y=0$. On the
other hand, the sources can be distributed over the hemisphere surface and
in the extension of the domain (see Fig. \ref{Fig2}). Figure \ref{Fig3}
shows the relative $L_{2}$-norm errors, $\epsilon _{2}$, 
\begin{equation}
\epsilon _{2}=\frac{\left\| f-f_{ref}\right\| _{2}}{\left\| f_{ref}\right\|
_{2}},  \label{co5.1}
\end{equation}%
where $f$ is the computed and $f_{ref}$ is the reference solution for all
computed source-receiver pairs as the filled contours in plane $\left[
R_{e}/R_{0},p\right] $. It is seen that the boundaries of the regions agree
well with the error approximation $R_{e}=R_{0}\epsilon _{2}^{-1/p}$. Figure %
\ref{Fig4} is computed by finding for each $R_{e}/R_{0}$ the smallest $p$ at
which the computed $\epsilon _{2}$ is smaller than the prescribed accuracy
and plotting the wall-clock time required for kernel computation via
factorization. It is seen that dependences computed at different prescribed $%
\epsilon _{2}$ have global minima in the range $2\leqslant
R_{e}/R_{0}\leqslant 3.5$ at $10^{-2}\geqslant \epsilon _{2}\geqslant 10^{-8}
$. This agrees well with the theoretical estimates provided above. Note that
non-monotonous character of the curves shown on this plot is due to discrete
change of parameters $R_{e}/R_{0}$ and $p$ and thresholding (the curves show
the prescribed error curves, not the actual errors, which can be close to or
several times smaller than the prescribed errors). Also, this figure shows
that optimization of the extended domain size can speed up computations
several times. 

\subsection{Optimization of the Boundary Element Solution}

The extension of the domain can substantially impact the complexity and
memory requirements of the BEM and it is preferable not to extend the domain
too much. Assume that $M=N$ and $R_{e}=R_{0}\left( 1+\delta \right) $, $%
\delta \ll 1$. In this case, $p\sim \delta ^{-1}\ln \left( 1/\epsilon\right) 
$, the first term in Eq. (\ref{co2}) can be neglected, and we obtain 
\begin{equation}
C_{fact}\sim B\frac{N}{\delta ^{3}}\ln ^{3}\frac{1}{\epsilon }\left( 1+\frac{%
\pi R_{0}^{2}}{A_{0}}\frac{2C}{B}\frac{\delta ^{2}}{\ln \frac{1}{\epsilon }}%
\right) \sim B\frac{N}{\delta ^{3}}\ln ^{3}\frac{1}{\epsilon }.  \label{co6}
\end{equation}
Here we assumed that the density of sampling of the domain extension is the
same as the density of sampling of surface $S\cup S_{g0}$ of area $A_{0}$,
means $\rho =N/A_{0}$. The last term in the parentheses then can be
neglected at $\delta \rightarrow 0$, due to $\pi R_{0}^{2}\leqslant A_{0}$, $%
\ln \frac{1}{\epsilon }>1,$ and assumed $C\sim B.$

There are different versions of the BEM including those using iterative
methods and the FMM, which have different computational costs, e.g., $%
O\left( N^{\alpha }\right) $, $\alpha \geqslant 1$, for problems of size $N.$
So, the cost of the BEM with the extended domain kernel computing at $\delta
\ll 1$ can be estimated as 
\begin{eqnarray}
C_{BEM} &=&C_{fact}+DN_{e}^{\alpha }=C_{fact}+DN^{\alpha }\left( 1+\frac{%
N_{e}-N}{N}\right) ^{\alpha }  \label{co7} \\
&\sim &C_{fact}+DN^{\alpha }\left( 1+\alpha \frac{\pi R_{0}^{2}}{A_{0}}%
\left( \frac{R_{e}^{2}}{R_{0}^{2}}-1\right) \right) \sim DN^{\alpha }+B\frac{%
N}{\delta ^{3}}\ln ^{3}\frac{1}{\epsilon }+2DN^{\alpha }\alpha \frac{\pi
R_{0}^{2}}{A_{0}}\delta .  \notag
\end{eqnarray}%
Here $D$ is some asymptotic constant related to the BEM cost. While it
appears that the cost of the last term is much smaller than the cost of the
first term, we retained it, as the first term is constant, while the last
term reflects the main mechanism of increase of the cost at the increasing $%
\delta $. Indeed, taking derivative of this expression we can see that $%
C_{BEM}$ has a global minimum, 
\begin{eqnarray}
\frac{dC_{BEM}}{d\delta } &=&-3B\frac{N}{\delta ^{4}}\ln ^{3}\frac{1}{%
\epsilon }+2\alpha D\frac{\pi R_{0}^{2}}{A_{0}}N^{\alpha },\quad \frac{%
dC_{BEM}}{d\delta }\left( \delta _{opt}\right) =0,\quad  \label{co9} \\
\delta _{opt} &=&\frac{3B}{2\alpha D}\frac{A_{0}}{\pi R_{0}^{2}}\left( \frac{%
1}{N^{\alpha -1}}\ln ^{3}\frac{1}{\epsilon }\right) ^{1/4}=O\left( \frac{1}{%
N^{\left( \alpha -1\right) /4}}\ln ^{3/4}\frac{1}{\epsilon }\right) .  \notag
\end{eqnarray}%
This shows that, indeed, at $\alpha >1$ and $N\rightarrow \infty $ and fixed
other parameters the optimal size of the extension, $\delta _{opt}$, turns
to zero. For the regular BEM with exact solvers we have $\alpha =3$.
However, when $N$ is not very high the cost of the FMM can be limited by the
cost of boundary integral computing, in which case $\alpha =2$. Also, this
can be achieved with iterative solvers. In the case of the FMM acceleration $%
\alpha =1$ or close to 1. In this case, $\delta _{opt}$ should not depend on 
$N$, and the only reason why $\delta _{opt}\ll 1$ can be related to large
values of asymptotic constant $D$. Usually, this is true for the FMM.

\section{BEM computing}

\begin{figure}[tbh]
\par
\begin{center}
\includegraphics[width=0.9\textwidth, trim=0.5in 1.95in 0.5in 0.5in]{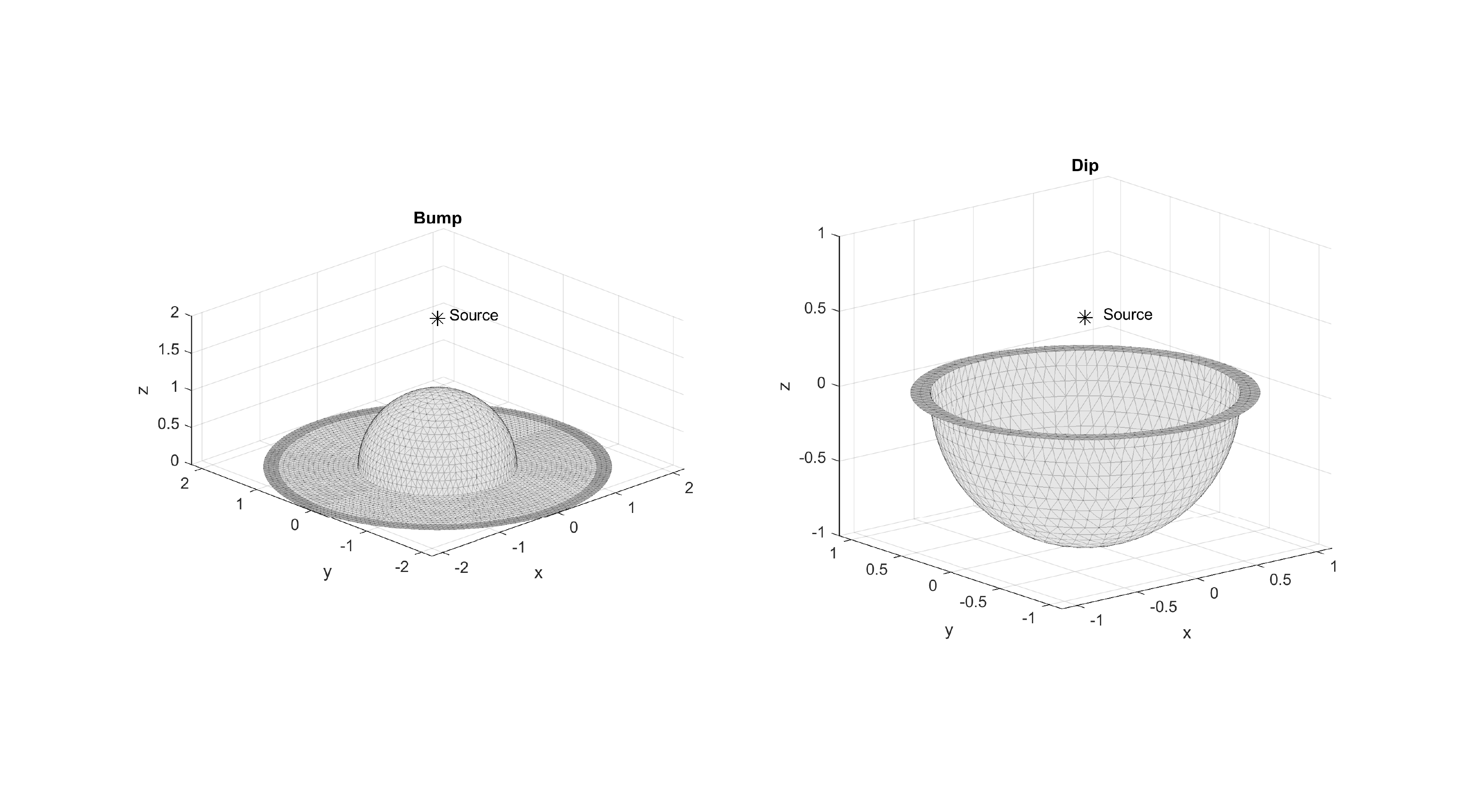}
\end{center}
\caption{The meshes of ``bump'' and ``dip'' and source locations used in
computations of Figs \ref{Fig7} and \ref{Fig8}. The parts corresponding to $%
S_{g0}$ are shaded in the light gray. The darker gray shade shows the parts
corresponding to $S_{ge}\backslash S_{g0}$. }
\label{Fig5}
\end{figure}

For illustrations of the kernel can be used with the boundary elements we
consider the following benchmark problems for validation. We assumed that
the ground elevation $z$ is the following function of $x$ and $y$, 
\begin{equation}
S_{g}:\left\{ (x,y,z),\quad z=\left\{ 
\begin{array}{c}
\pm \sqrt{1-\rho ^{2}},\quad \rho <1, \\ 
0,\quad \rho \geqslant 1.%
\end{array}%
\right. ,\quad \rho =\sqrt{x^{2}+y^{2}}\right\} .  \label{e7}
\end{equation}%
Here the sign ``+'' corresponds to a ``bump'', while the sign ``-''
corresponds to a ``dip''. Assume then that the electric field is generated
by a monopole source (charge) of unit intensity located at $\mathbf{x}%
_{s}=\left( 0,0,h\right) $, and the domain $\Omega _{0}$ is bounded by a
sphere of radius $R_{0}=h$ ($h>1$) centered at the origin for the ``bump''
and , $R_{0}=1$ for the ``dip'' ($\left| h\right| <1$). Figure \ref{Fig5}
shows the meshes used in computations.

For validation purposes, the solution of the problem for the bump can be
found analytically using the method of images. In this case, the image of
the bump is a hemisphere, which together with the ``image'' bump forms a
sphere of radius 1 on which the potential is zero. The image source of
negative intensity (sink) is located at $\mathbf{x}_{s}^{\ast }=\left(
0,0,-h\right) $, and we have the following expression for the potential 
\begin{equation}
\phi \left( \mathbf{y}\right) =G\left( \mathbf{y},\mathbf{x}_{s}\right)
-G\left( \mathbf{y},\mathbf{x}_{s}^{\ast }\right) -\frac{1}{h}G\left(\mathbf{%
y},\mathbf{x}_{r}^{+}\right) +\frac{1}{h}G\left( \mathbf{y},\mathbf{x}%
_{r}^{-}\right) ,\quad \mathbf{x}_{r}^{\pm }=\left( 0,0,\pm \frac{1}{h}%
\right) ,  \label{e8}
\end{equation}
where $\mathbf{x}_{r}^{+}$ and $\mathbf{x}_{r}^{-}$ are the images of $%
\mathbf{x}_{s}$ and $\mathbf{x}_{s}^{\ast }$ with respect to the sphere. We
are not aware about analytical solutions for the case of ``dip'' for which
the method of images cannot be applied.

The boundary integral solution of the problem can be represented in the form 
\begin{equation}
\phi \left( \mathbf{y}\right) =G\left( \mathbf{y},\mathbf{x}_{s}\right)
+K^{(D)}\left( \mathbf{y},\mathbf{x}_{s}\mathbf{;}R_{e}\right)+\int_{S_{ge}}%
\sigma ^{(D)}\left( \mathbf{x}\right) \left[ G\left( \mathbf{y}-\mathbf{x}%
\right) +K^{(D)}\left( \mathbf{y},\mathbf{x;}R_{e}\right) \right] dS\left( 
\mathbf{x}\right) ,  \label{e9}
\end{equation}
where the charge density can be found from 
\begin{equation}
\int_{S_{ge}}\sigma ^{(D)}\left( \mathbf{x}\right) \left[ G\left( \mathbf{y},%
\mathbf{x}\right) +K^{(D)}\left( \mathbf{y},\mathbf{x;}R_{e}\right) \right]%
dS\left( \mathbf{x}\right) =-G\left( \mathbf{y},\mathbf{x}_{s}\right)
-K^{(D)}\left( \mathbf{y},\mathbf{x}_{s}\mathbf{;}R_{e}\right) ,\quad 
\mathbf{y}\in S_{ge}.  \label{e10}
\end{equation}

In the case of ``bump'' we computed and compared four solutions: 1)
analytical; 2) regular BEM using method of images (two symmetric sources and
a sphere) (``BEMimage''); 3) regular BEM, where the domain ground plane was
simply truncated and the effect of the infinite plane ignored (``BEM''); 4)
BEM with the kernel accounting for the infinite plane as described by Eqs (%
\ref{e9}) and (\ref{e10}) (``BEMinf''). In case of ``dip'' only solutions 3)
and 4) were computed. In all cases we varied $R_{e}$. The truncation number $%
p$ for ``BEMinf'' method was selected based on the theoretical estimate, $%
p=\log \epsilon /\log (R_{0}/R_{e})$, where the prescribed error $\epsilon $
varied in range $10^{-3}-10^{-6}$.

\begin{figure}[tbh]
\vspace{-4pt}
\par
\begin{center}
\includegraphics[width=0.9\textwidth, trim=0.5in 0.85in 0.5in 0.in]{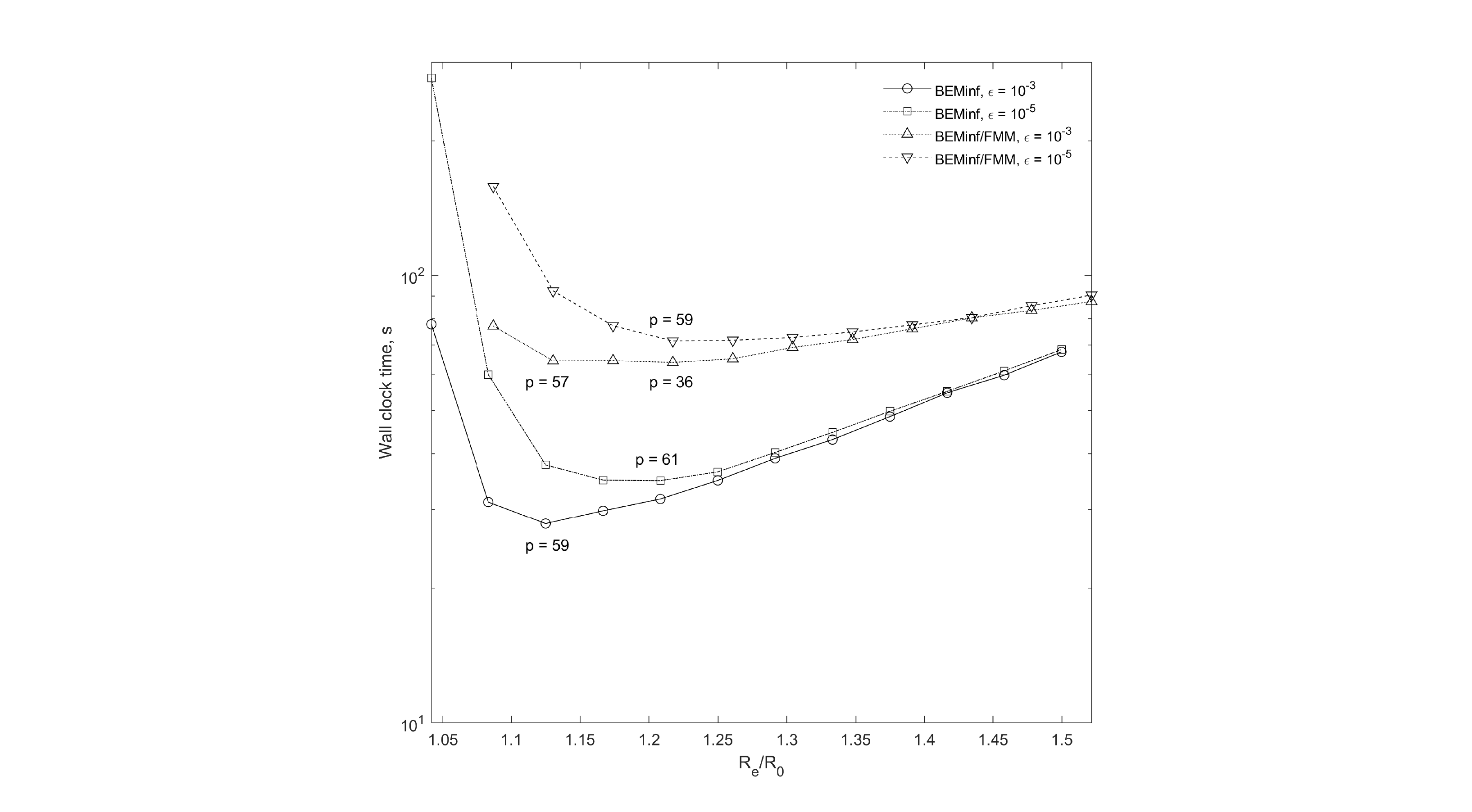}
\end{center}
\caption{The wall-clock times (in seconds) for solution of the ``bump''
problem shown in Fig. \ref{Fig5} for different $R_{e}/R_{0}$, different
prescribed accuracies, and different versions of the BEMinf (with and
without the FMM acceleration). The values of $p$ corresponding to the minima
of the dependences are shown near the curves. In the regular BEMinf, the
number of faces for surface $S_{g0}$ is $N_{0}=3594$ and varied in the range 
$N_{e}=3901-8106$ for surface $S_{ge}.$ In the BEMinf accelerated with the
FMM the numbers are $N_{0}=13244$ and $N_{e}=16568-30717$, respectively. 
\newline
}
\label{Fig6}
\end{figure}
\begin{figure}[tbh]
\vspace{-4pt}
\par
\begin{center}
\includegraphics[width=0.9\textwidth, trim=0.5in 0.75in 0.5in 0.in]{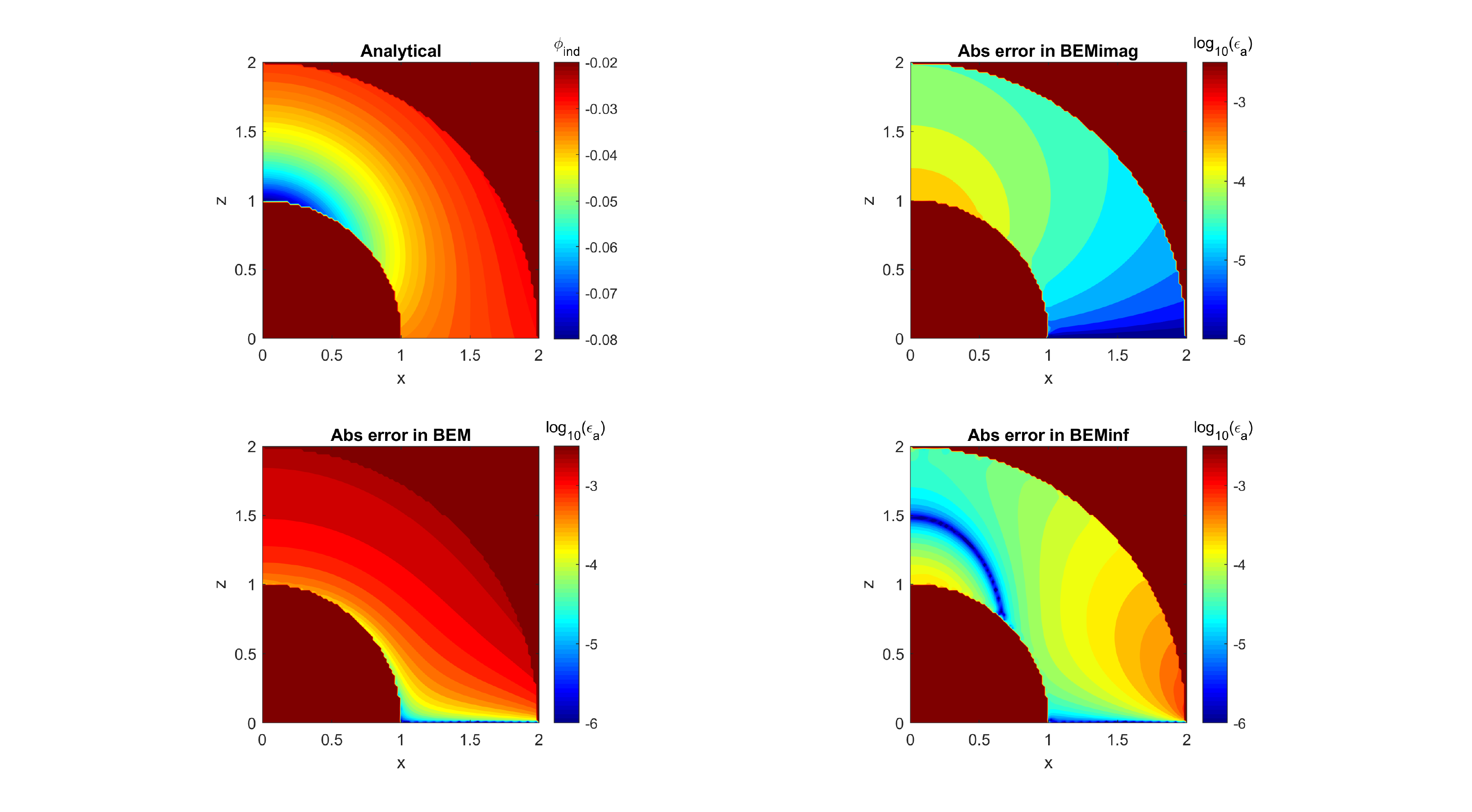}
\end{center}
\caption{The induced electric potential, $\protect\phi _{ind}\left( \mathbf{y%
}\right) =\protect\phi \left( \mathbf{y}\right) -G\left( \mathbf{y},\mathbf{x%
}_{s}\right) $, according to the analytical solution (\ref{e8}) , and the
absolute errors in domain $\Omega _{0}$ obtained using different methods of
solution for the ``bump'' case, $h=2$, $\protect\delta =0.0935$, $N_{0}=6401$%
, and $N_{e}=7661$. The mesh is shown in Fig. \ref{Fig5} (on the left).}
\label{Fig7}
\end{figure}

Figure \ref{Fig6} illustrates the effect of the domain extension on the
performance for different prescribed accuracies and the BEM complexities for
the ``bump'' case ($h=2$). Note that here the FMM accelerated BEM was
implemented in the way that the boundary integrals were computed within some
cutoff radius of the evaluation point using Matlab, while the FMM procedure
was executed in the form of a high performance library, so the cost of the
FMM BEM is mainly determined by the Matlab integral computations, which is
consistent with the kernel evaluation also implemented in Matlab. Also, in
the regular BEM the solution time was much smaller than the integral
computation time. So, the methods shown can be characterized as $O\left(
N^{2}\right) $ for the BEM and $O\left( N\right) $ for the BEM/FMM. It is
seen that the curves form minima (sometimes a few local minima), at
relatively low $\delta =R_{e}/R_{0}-1$. According to the theory (see Eq. (%
\ref{co9})), increasing the prescribed accuracy results in an increase of
the optimal $R_{e}$, which is seen on the figure. For $\delta $
substantially larger than $\delta _{opt}$ the cost of kernel computation
becomes negligible in comparison to the BEM costs, and the curves computed
for different $\epsilon $ become close to each other.

Figure \ref{Fig7} shows the analytical solution (\ref{e8}) (the induced
part, $\phi _{ind}\left( \mathbf{y}\right) =\phi \left( \mathbf{y}\right)
-G\left( \mathbf{y},\mathbf{x}_{s}\right) $), and the absolute errors in the
computational domain obtained using different methods of solution for the
``bump'' case, $h=2$, $\delta =0.0935$, $N_{0}=6401$, and $N_{e}=7661$,
where $N_{0}$ and $N_{e}$ are the number of boundary elements in $\Omega
_{0} $ and $\Omega _{e}$, respectively. The numerical solution using image
sources is the most accurate and its relative $L_{2}$-norm error is $%
\epsilon _{2}=2.3\cdot 10^{-3}$. The BEMinf shows error $\epsilon
_{2}=4.5\cdot 10^{-3}$ (at the prescribed accuracy $\epsilon =10^{-4}$)
while the BEM shows $\epsilon _{2}=3.7\cdot 10^{-2}$, which is one order of
magnitude larger. This figure shows also different error distributions for
different methods inside the computational domain.

\begin{figure}[tbh]
\vspace{-4pt}
\par
\begin{center}
\includegraphics[width=0.9\textwidth, trim=0.5in 1.75in 0.5in 0.in]{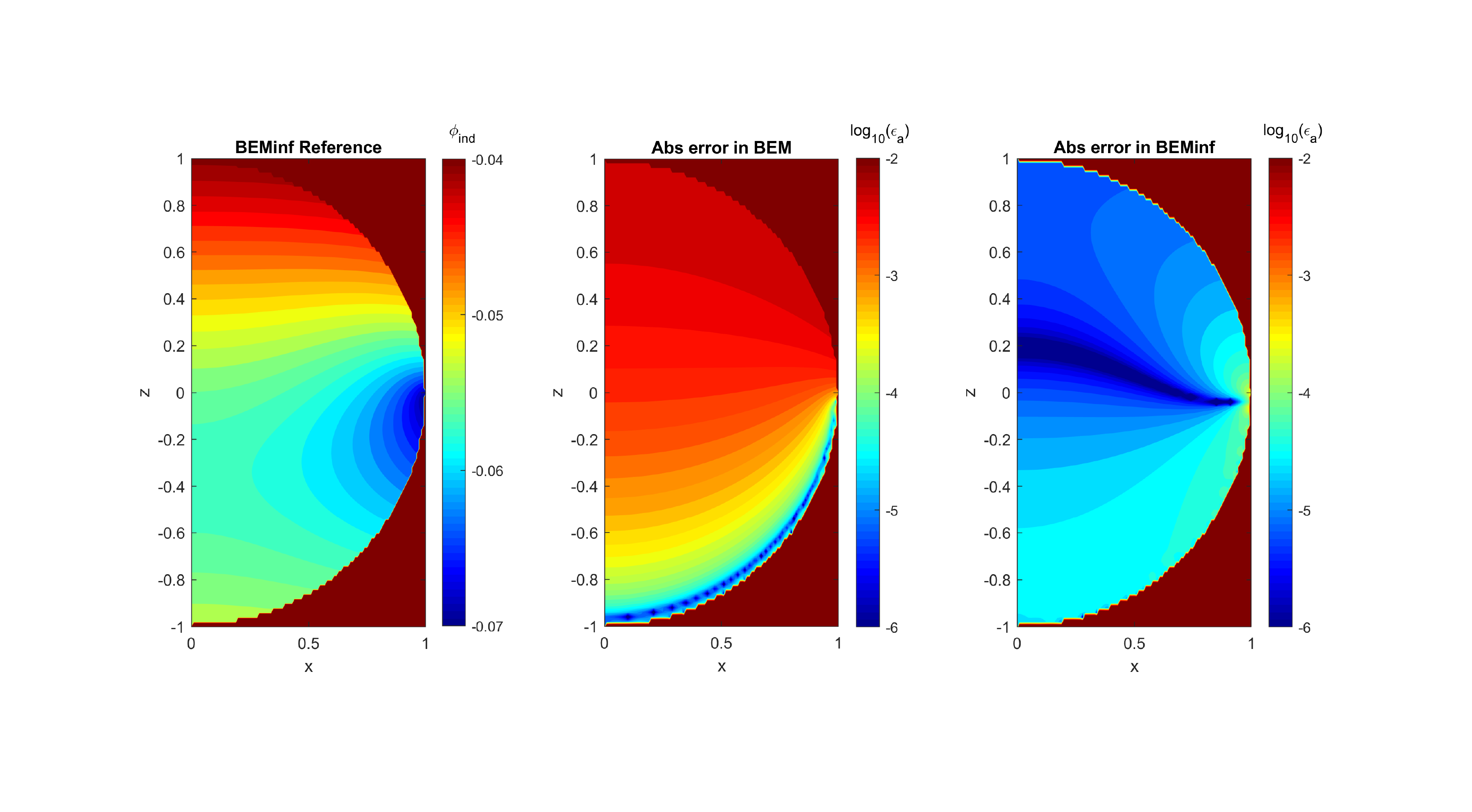}
\end{center}
\caption{The induced electric potential, $\protect\phi _{ind}\left( \mathbf{y%
}\right) =\protect\phi \left( \mathbf{y}\right) -G\left( \mathbf{y},\mathbf{x%
}_{s}\right) $, for the reference solution computed using BEMinf with
parameters $R_{e}/R_{0}=1.5$, $N_{0}=6401$, $N_{e}=14038,$ and prescribed $%
\protect\epsilon =10^{-6}$, and the absolute errors in domain $\Omega _{0}$
obtained using different methods of solution for the ``dip'' case, $h=0.5$, $%
\protect\delta =0.124$, $N_{0}=1592$, and $N_{e}=2017$. The mesh is shown in
Fig. \ref{Fig5} (on the right).}
\label{Fig8}
\end{figure}
\begin{figure}[tbh]
\vspace{-4pt}
\par
\begin{center}
\includegraphics[width=0.9\textwidth, trim=0.5in 0.75in 0.5in 0.in]{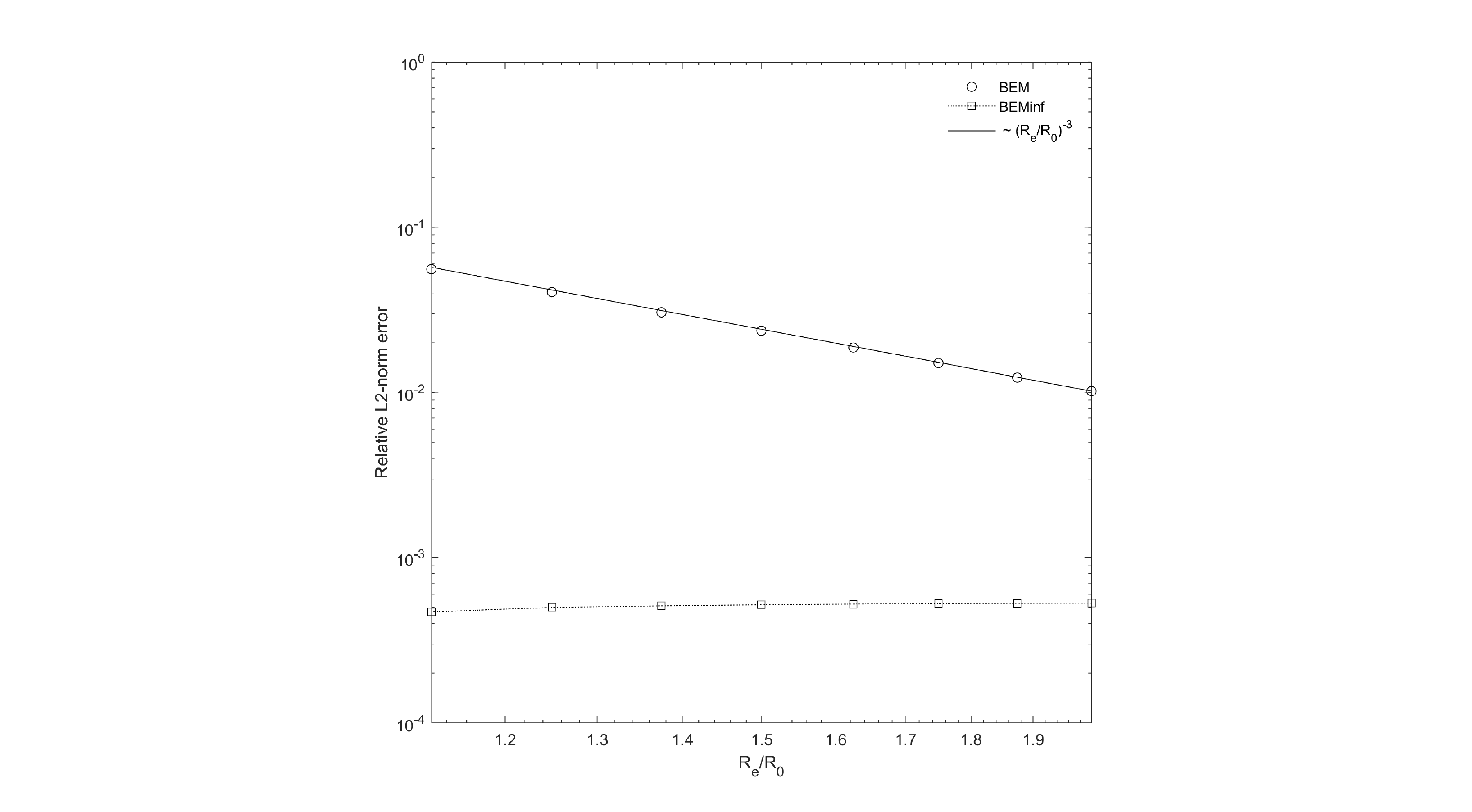}
\end{center}
\caption{The relative $L_{2}$-norm error $\protect\epsilon _{2}$ for the
``dip'' case as a function of parameter $R_{e}/R_{0}$ for BEMinf and regular
BEM. The solid line shows the dependence $\protect\epsilon _{2}=C\left(
R_{e}/R_{0}\right) ^{-3}$ with the constant $C$ fiitting the error of
regular BEM at $R_{e}/R_{0}=2$. }
\label{Fig9}
\end{figure}
Figure \ref{Fig8} shows the reference solution $\left( \phi _{ind}\left( 
\mathbf{y}\right) \right) $ for the ``dip'' case and the absolute errors in
the computational domain obtained using BEMinf and the regular BEM, $h=0.5$, 
$\delta =0.124$, $N_{0}=1592$, and $N_{e}=2017$. The reference solution was
computed using BEMinf with parameters $\delta =0.5$, $N_{0}=6401$, and $%
N_{e}=14038$ and prescribed $\epsilon =10^{-6}$. The BEMinf shows error $%
\epsilon _{2}=4.7\cdot 10^{-4}$ (at the prescribed accuracy $\epsilon
=10^{-4}$) while the BEM shows $\epsilon _{2}=5.6\cdot 10^{-2}$, which is
two orders of magnitude larger. Similarly to the case of ``bump'' one can
observe different error distributions for the BEMinf and the BEM. Similarly
to the case of ``bump'' the largest BEMinf errors are observed near the edge
of domain $\Omega _{0}$ coinciding with the ground plane.

Figure \ref{Fig9} plots $\epsilon _{2}$ for the ``dip'' case as a function
of parameter $R_{e}/R_{0}$. It is seen that the error is practically
constant for BEMinf, while it decays inversely proportionally to the volume
of the computational domain for the regular BEM. It is clear then that the
results of BEMinf can be obtained by simple extension of the computational
domain, but for substantially higher computational costs. For example, to
achieve two orders of magnitude in error reduction one should have $%
R_{e}\sim 5R_{0}$, which may require one order of magnitude increase in the
number of boundary elements and increase the computational costs 10-1000
times, based on the solution method used (scaled as $O\left( N^{\alpha
}\right) $, $1\leqslant \alpha \leqslant 3$).

\section{Conclusion}

We introduced analogues of Green's function for an infinite plane with a
circular hole. We also proposed and tested efficient methods for computing
of these functions, which may have broad applications to model infinite
ground or ocean surfaces in the problems requiring Laplace equation solvers.
While such functions can be expressed in the forms of finite integrals,
their computation can be costly. However, a computationally efficient
factorization of these functions can be done, which drastically reduces the
computational costs, and, in fact, creates rather small overheads to the
problems, which simply ignore the existence of ground outside the
computational domain. Also, we conducted some optimization study and
provided actual boundary element computations for some benchmark problems.
The study shows that the accuracy of the solution using the analogues of
Green's function improves substantially compared to a simple domain
truncation, which also can become computationally costly for improved
accuracy computations.

\section{Acknowledgements}

This work is supported by Cooperative Research Agreement (W911NF1420118)
between the University of Maryland and the Army Research Laboratory, with
David Hull, Ross Adelman and Steven Vinci as Technical monitors.

\appendix

\section{Derivation of Recurrence Relations}

Recurrence (\ref{ed4.3}) can be derived by comparing two representations of
integral%
\begin{eqnarray}
I_{n}^{m} &=&\frac{1}{\xi ^{n+1}}\int_{0}^{2\pi }e^{im\varphi }d\varphi
\int_{0}^{\xi }\left( 1-2\zeta \cos \varphi +\zeta ^{2}\right) ^{1/2}\zeta
^{n}d\zeta  \label{A1} \\
&=&u_{n}^{m}+\xi ^{2}u_{n+2}^{m}-\xi
\left(u_{n+1}^{m-1}+u_{n+1}^{m+1}\right) ,  \notag
\end{eqnarray}%
which can be checked using definition of $u_{n}^{m}\left( \xi \right) $ (\ref%
{ed3}), and integration by parts: 
\begin{eqnarray}
I_{n}^{m} &=&\frac{1}{\left( n+1\right) \xi ^{n+1}}\int_{0}^{2\pi}e^{im%
\varphi }d\varphi \int_{0}^{\xi }\left( 1-2\zeta \cos \varphi
+\zeta^{2}\right) ^{1/2}d\zeta ^{n+1}  \label{A2} \\
&=&\frac{1}{n+1}\left[ v_{m}-\xi ^{2}u_{n+2}^{m}+\frac{1}{2}\xi
\left(u_{n+1}^{m-1}+u_{n+1}^{m+1}\right) \right] ,\quad
v_{m}=\int_{0}^{2\pi}e^{im\varphi }\left( 1-2\xi \cos \varphi +\xi
^{2}\right) ^{1/2}d\varphi .  \notag
\end{eqnarray}
Equation (\ref{ed4.4}) can be checked directly using integral representation
of $v_{m}$ and definition of $w_{m}\left( \xi \right) .$

To obtain recurrence (\ref{ed4.5}) we note the following transformation 
\begin{gather}
i\int_{0}^{2\pi }\frac{e^{i\left( m-1\right) \varphi }\sin \varphi d\varphi 
}{\left( 1-2\xi \cos \varphi +\xi ^{2}\right) ^{1/2}}=\frac{i}{\xi }%
\int_{0}^{2\pi }e^{i\left( m-1\right) \varphi }d\left( 1-2\xi \cos \varphi
+\xi ^{2}\right) ^{1/2}  \label{A3} \\
=-\frac{i}{\xi }\int_{0}^{2\pi }\left( 1-2\xi \cos \varphi +\xi
^{2}\right)^{1/2}de^{i\left( m-1\right) \varphi }=\left( m-1\right) \frac{1}{%
\xi }\int_{0}^{2\pi }\frac{e^{i\left( m-1\right) \varphi }\left( 1-2\xi
\cos\varphi +\xi ^{2}\right) }{\left( 1-2\xi \cos \varphi +\xi ^{2}\right)
^{1/2}}d\varphi .  \notag
\end{gather}
Using this expression, we obtain 
\begin{eqnarray}
w_{m} &=&\int_{0}^{2\pi }\frac{e^{i\left( m-1\right) \varphi
}e^{i\varphi}d\varphi }{\left( 1-2\xi \cos \varphi +\xi ^{2}\right) ^{1/2}}%
=\int_{0}^{2\pi }\frac{e^{i\left( m-1\right) \varphi }\left( \cos
\varphi+i\sin \varphi \right) d\varphi }{\left( 1-2\xi \cos \varphi
+\xi^{2}\right) ^{1/2}}  \label{A4} \\
&=&\left( m-1\right) \frac{1+\xi ^{2}}{\xi }\int_{0}^{2\pi }\frac{%
e^{i\left(m-1\right) \varphi }d\varphi }{\left( 1-2\xi \cos \varphi +\xi
^{2}\right)^{1/2}}-\frac{2m-3}{2}\int_{0}^{2\pi }\frac{\left(
e^{i\varphi}+e^{-i\varphi }\right) e^{i\left( m-1\right) \varphi }}{\left(
1-2\xi \cos\varphi +\xi ^{2}\right) ^{1/2}}d\varphi  \notag \\
&=&\left( m-1\right) \frac{1+\xi ^{2}}{\xi }w_{m-1}-\frac{2m-3}{2}%
\left(w_{m}+w_{m-2}\right) .  \notag
\end{eqnarray}
It is not difficult to see that Eq. (\ref{ed4.5}) follows from this
expression.

Functions $w_{0}$ and $w_{1}$ can be expressed immediately via complete
elliptic integrals. Indeed, we have 
\begin{eqnarray}
w_{m}\left( \xi \right) &=&\int_{\pi }^{3\pi }\frac{e^{im\varphi }d\varphi }{%
\left( 1-2\xi \cos \varphi +\xi ^{2}\right) ^{1/2}}=\left( -1\right)
^{m}\int_{0}^{2\pi }\frac{e^{im\varphi }d\varphi }{\left( 1+2\xi \cos
\varphi +\xi ^{2}\right) ^{1/2}}  \label{A5} \\
&=&2\left( -1\right) ^{m}\int_{0}^{\pi }\frac{\cos m\varphi d\varphi }{%
\left( 1+2\xi \cos \varphi +\xi ^{2}\right) ^{1/2}}=\frac{4\left( -1\right)
^{m}}{1+\xi }\int_{0}^{\pi /2}\frac{\cos \left( 2m\psi \right) d\psi }{%
\left( 1-\mu \sin ^{2}\psi \right) ^{1/2}},\quad \mu =\frac{4\xi }{\left(
1+\xi \right) ^{2}}.  \notag
\end{eqnarray}%
So, 
\begin{equation}
w_{0}\left( \xi \right) =\frac{4}{1+\xi }K\left( \mu \right) ,\quad
w_{1}\left( \xi \right) =\frac{4}{1+\xi }\left[ \left( \frac{2}{\mu }%
-1\right) K\left( \mu \right) -\frac{2}{\mu }E\left( \mu \right) \right] ,
\label{A6}
\end{equation}%
where $K$ and $E$ are defined by Eq. (\ref{ed4.7}). The Landen
transformation of elliptic functions \cite{Abramowitz1964book} shows that 
\begin{eqnarray}
K\left( \mu \right) &=&\frac{2}{1+\mu _{1}^{1/2}}K\left( \mu _{2}\right)
,\quad \mu _{1}=1-\mu ,\quad \mu _{2}=\left( \frac{1-\mu _{1}^{1/2}}{1+\mu
_{1}^{1/2}}\right) ^{2}  \label{A7} \\
\quad E\left( \mu \right) &=&\left( 1+\mu _{1}^{1/2}\right) E\left( \mu
_{2}\right) -\frac{2\mu _{1}^{1/2}}{1+\mu _{1}^{1/2}}K\left( \mu _{2}\right)
.  \notag
\end{eqnarray}%
Substituting this into Eq. (\ref{A6}) and using the value of $\mu $ from Eq.
(\ref{A5}) we obtain Eq. (\ref{ed4.6}).

Note now that relation (\ref{ed4.3}) is also valid in case $m=0$. In this
case we can use symmetry $u_{n}^{-1}=u_{n}^{1}$. Since Eq. (\ref{ed4.4})
also holds in this case, we obtain using symmetry $w_{-1}=w_{1}$ and Eq. (%
\ref{ed4.6}), 
\begin{equation}
v_{0}=\left( 1+\xi ^{2}\right) w_{0}-2\xi w_{1}=8E\left( \xi ^{2}\right)
-4\left( 1-\xi ^{2}\right) K\left( \xi ^{2}\right) ,  \label{A8}
\end{equation}
and confirms recurrence (\ref{ed4.2}).

So, the last relation to prove is Eq. (\ref{ed4.1}). This relation can be
obtained using differentiation relations for the complete elliptic integral, 
\begin{equation}
K\left( \mu \right) =\frac{d}{d\mu }\left[ \mu \left( 1-\mu \right) \frac{%
dK\left( \mu \right) }{d\mu }\right] ,\quad \frac{dK\left( \mu \right) }{%
d\mu }=\frac{E\left( \mu \right) }{2\mu \left( 1-\mu \right) }-\frac{%
K\left(\mu \right) }{2\mu },  \label{A9}
\end{equation}
and integration by parts, 
\begin{eqnarray}
u_{n}^{0} &=&\frac{4}{\xi ^{n+1}}\int_{0}^{\xi }\zeta ^{n}K\left(
\zeta^{2}\right) d\zeta =\frac{8}{\xi ^{n+1}}\int_{0}^{\xi ^{2}}\mu
^{\left(n-1\right) /2}\frac{d}{d\mu }\left( \mu \left( 1-\mu \right) \frac{%
dK\left(\mu \right) }{d\mu }\right) d\mu  \label{A10} \\
&=&\frac{8}{\xi ^{n+1}}\left[ \xi ^{n+1}\left( 1-\xi ^{2}\right) \frac{%
dK\left( \xi ^{2}\right) }{d\left( \xi ^{2}\right) }-\frac{n-1}{2}%
\int_{0}^{\xi ^{2}}\mu ^{\left( n-1\right) /2}\left( 1-\mu \right)
dK\left(\mu \right) \right]  \notag \\
&=&\frac{4}{\xi ^{2}}\left[ E\left( \xi ^{2}\right) -n\left(
1-\xi^{2}\right) K\left( \xi ^{2}\right) +\frac{\left( n-1\right) ^{2}}{4}%
u_{n-2}^{0}\right] -\left( n^{2}-1\right) u_{n}^{0}.  \notag
\end{eqnarray}
Equation (\ref{ed4.1}) is just a rearrangement of terms in this relation.


\section{References}

\begin{thebibliography}{99}
\bibitem{Adelman2017:SC} R. Adelman, Extremely large, wide-area power-line
models, Supercomputing 2017,
https://sc17.supercomputing.org/SC17\%20Archive/tech\_poster/poster%
\_files/post251s2-file3.pdf.

\bibitem{DAmore96:IEEE} M. D'Amore, M. S. Sarto, Simulation models of a
dissipative transmission line above a lossy ground for a wide-frequency
range. I. Single conductor configuration, IEEE Trans. Electromag. Compat.
38(2) (1996) 127-138. https://doi.org/10.1109/15.494615.

\bibitem{Trlep2009:IEEE} M. Trlep, A. Hamler, M. Jesenik, B. Stumberger,
Electric field distribution under transmission lines dependent on ground
surface, IEEE Trans. Mag. 45(3) (2009) 1748-1751.
https://doi.org/10.1109/TMAG.2009.2012806.

\bibitem{Zhang2012:Electrostat} Q. Zhang, J. Yang, D. Li, Z. Wang,
Propagation effects of a fractal rough ocean surface on the vertical
electric field generated by lightning return strokes, J. Electrostat.
70(1)(2012) 54-59. https://doi.org/10.1016/j.elstat.2011.10.003.

\bibitem{Chen2003:IEEE} K. S. Chen, T.-D. Wu, L. Tsang, Q. Li, J. Shi, A.K.
Fung, Emission of rough surfaces calculated by the integral equation method
with comparison to three-dimensional moment method simulations, IEEE Trans.
Geo. Rem. Sensing 41(1) (2003) 90-101.
https://doi.org/10.1109/TGRS.2002.807587.

\bibitem{Givoli1992:book} D. Givoli, Numerical Methods for Problems in
Infinite Domains, Elsevier, Amsterdam, 1992.

\bibitem{Tsinkov1998:ANM} S.V. Tsynkov, Numerical solution of problems on
unbounded domains. A review, Appl. Num. Math. 27(4) (1998) 465-532.
https://doi.org/10.1016/S0168-9274(98)00025-7.

\bibitem{Jackson1998:Wiley} J.D. Jackson, Classical Electrodynamics, third
ed., John Wiley and Sons, New York, 1998.

\bibitem{Crowdy2007:IMA} D. Crowdy, J. Marshall, Green's functions for
Laplace's equation in multiply connected domains, IMA J. Appl. Math. 72
(2007) 278-301. https://doi.org/10.1093/imamat/hxm007.

\bibitem{Gumerov2014JCP} N.A. Gumerov and R. Duraiswami, A method to compute
periodic sums, J. Comput. Phys. 272 (2014) 307-326.
https://doi.org/10.1016/j.jcp.2014.04.039.

\bibitem{Moroz2006:JP} A. Moroz, Quasi-periodic Green's functions of the
Helmholtz and Laplace equations, J. Phys. A: Math. Gen. 39(36) (2006) 11247.
https://doi.org/10.1088\%2F0305-4470\%2F39\%2F36\%2F009.

\bibitem{Bao2013:IPI} G. Bao, J. Lin, Near-field imaging of the surface
displacement on an infinite ground plane, Inverse Problems \& Imaging, 7(2)
(2013), 377-396. https:/doi.org/10.3934/ipi.2013.7.377.

\bibitem{Chandler-Wilde2006:SIAM} S. N. Chandler-Wilde, E. Heinemeyer, R.
Potthast, Acoustic scattering by mildly rough unbounded surfaces in three
dimensions, SIAM J. Appl. Math. 66 (2006) 1002--1026.
https://doi.org/10.1137/050635262.

\bibitem{Holmes2015:MTSJ} J.J. Holmes, Past, present, and future of
underwater sensor arrays to measure the electromagnetic field signatures of
naval vessels, Marine Tech. Soc. J. 49(6) (2015) 123-133.
https://doi.org/10.4031/MTSJ.49.6.1.

\bibitem{Yue2016:IEEE} R. Yue, P. Hu, J. Zhang, The influence of the
seawater and seabed interface on the underwater low frequency
electromagnetic field signatures, 2016 IEEE/OES China Ocean Acoustics (COA),
Harbin (2016) 1-7. https://doi.org/10.1109/COA.2016.7535694.

\bibitem{Wang2018:PIER} X. Wang, Q. Xu, J. Zhang, Simulating underwater
electric field signal of ship using the boundary element method, Progress In
Electromag. Res. M 76 (2018) 43-54. https:/doi.org/10.2528/PIERM18092706.

\bibitem{Rozhdestvensky2000:book} K.V. Rozhdestvensky, Aerodynamics of a
Lifting System in Extreme Ground Effect, Springer, Berlin-Heidelberg, 2000.

\bibitem{Duffy2015book} D. G. Duffy, Green's Functions with Applications,
CRC Press, Taylor \& Francis Group, FL, 2015.

\bibitem{Morse1953book} P.M. Morse, H. Feshbach, Methods of Theoretical
Physics, McGraw-Hill, NY, 1953.

\bibitem{Gumerov2005report} N.A. Gumerov, R. Duraiswamin, Comparison of the
efficiency of translation operators used in the fast multipole method for
the 3D Laplace equation, UMIACS TR 2005-09, Also issued as Computer Science
Technical Report CS-TR-\# 4701, University of Maryland, College Park, 2005.
https://drum.lib.umd.edu/bitstream/handle/1903/3023/LaplaceTranslation\_CSTR%
\_4701.pdf.

\bibitem{Abramowitz1964book} M. Abramowitz, I.A. Stegun, Handbook of
Mathematical Functions, National Bureau of Standards, Washington D.C., 1972.

\bibitem{Adelman2017IEEE} R. Adelman, N.A Gumerov, and R. Duraiswami,
FMM/GPU-accelerated boundary element method for computational magnetics and
electrostatics, IEEE Trans. Mag. 53(12) (2017), 7002311.
https://doi.org/10.1109/TMAG.2017.2725951.

\bibitem{Gumerov2008JCP} N.A. Gumerov and R. Duraiswami, Fast multipole
methods on graphics processors, J. Comput. Phys. 227 (2008), 8290-8313.
https://doi.org/10.1016/j.jcp.2008.05.023.
\end{thebibliography}
\end{document}